\documentclass[a4paper,11pt,reqno]{article}
\usepackage[T1]{fontenc}
\usepackage{courier}
\usepackage[scaled=0.92]{helvet}

\usepackage{amscd,amsmath,amssymb,amsthm,amsfonts,epsfig,graphics}

\usepackage[english]{babel}
\usepackage{amsthm,amsfonts,amssymb,amsmath}
\usepackage{hyperref}
\hypersetup{
colorlinks=true,
citecolor=red,
linkbordercolor={1 1 1},
citebordercolor={1 1 1},
}
\usepackage{cleveref}
\usepackage{algpseudocode,algorithm,float}
\usepackage{mathrsfs,caption}
\usepackage[percent]{overpic} 				
\usepackage[
placement=left,
angle=90,
color=black!60,
firstpage=false,
scale=2.5,
hshift=0,
vshift=110
]{background}
\backgroundsetup{contents={\tiny{}}}

\usepackage{xcolor} 

\theoremstyle{plain}
\newtheorem{satz}{Theorem}[section]
\newtheorem{prop}[satz]{Proposition}

\newtheorem{lem}[satz]{Lemma}

\theoremstyle{definition}

\newtheorem{rem}[satz]{Remark}

\newcommand{\mx}{\mbox}
\newcommand{\rw}{\rightarrow}
\newcommand{\de}{\displaystyle}
\newcommand{\ml}{\mathcal}

\newcommand{\pl}{\partial}

\newcommand{\beq}[1]{\begin{equation} \label{#1}}
\newcommand{\eeq}{\end{equation}}
\newcommand{\beqar}[1]{\beq{#1} \begin{array}{rcl}}
\newcommand{\eeqar}{\end{array} \eeq}

\newcommand{\drx}[1]{\check{#1}}
\newcommand{\ddrx}[1]{\check{\check{#1}}}

\newcommand{\ue}[1]{\underline{#1}}

\providecommand{\ep}{\varepsilon}

\providecommand{\RR}{\mathbb{R}}

\providecommand{\ZZ}{\mathbb{Z}}
\providecommand{\NN}{\mathbb{N}}

\providecommand{\TT}{\mathbb{T}}

\newcommand{\norm}[1]{\left \lVert#1 \right\rVert_{\infty}}
\newcommand{\enorm}[1]{\left \lVert#1 \right\rVert_{2}}

\newcommand{\K}[1]{\mathbf{K}({#1})}
\newcommand{\Kp}[1]{\mathbf{K}'({#1})}
\newcommand{\E}[1]{\mathbf{E}({#1})}

\newcommand{\Heq}[2]{\overset{\left(#1\right)}{\underset{}{#2}}}

\usepackage[utf8]{inputenc}

\hypersetup{pdfauthor={Fortunati, A.},%
            pdftitle={On the constructivity...},%
            pdfsubject={},%
            pdfkeywords={}%
}

\DeclareMathOperator{\dist}{dist}

\DeclareMathOperator{\diam}{diam}
\DeclareMathOperator{\diag}{diag}

\DeclareMathOperator{\Res}{Res}

\DeclareMathOperator{\cn}{cn}
\DeclareMathOperator{\dn}{dn}
\DeclareMathOperator{\sn}{sn}
\DeclareMathOperator{\am}{am}
\DeclareMathOperator{\sd}{sd}
\DeclareMathOperator{\sech}{sech}
\DeclareMathOperator{\Mod}{\mathrm{mod}}
\DeclareMathOperator{\trid}{tridiag}

\DeclareMathOperator{\mcd}{\mathrm{mcd}}
\DeclareMathOperator{\meas}{meas}

\makeatletter
\renewcommand*{\@fnsymbol}[1]{\ensuremath{\ifcase#1\or *\or \mathsection \or (b)\or \else \fi}}
\makeatother

\title{\LARGE{\textbf{On the constructivity of the variational approach to Arnold's Diffusion.}}}
\date{}
\author{%
Alessandro Fortunati\thanks{Current affiliation: Department of Electrical Engineering and Information Technology (DIETI), University of Naples Federico II, Italy. Email: alessandro.fortunati@unina.it.}}

\usepackage{graphicx}
\begin{document}

\maketitle

\begin{abstract}
The aim of this paper is to discuss the constructivity of the method originally introduced by U. Bessi to approach the phenomenon of topological instability commonly known as Arnold's Diffusion. By adapting results and proofs from existing works and introducing additional tools where necessary, it is shown how, at least for a (well known) paradigmatic model, it is possible to obtain a rigorous proof on a suitable discrete space, which can be fully implemented on a computer. A selection of explicitly constructed diffusing trajectories for the system at hand is presented in the final section.
 
\end{abstract}

\section{Background}\label{sec:background}
In his pioneering paper \cite{arn}, by exhibiting a cleverly chosen model of nearly-integrable Hamiltonian system with three degrees of freedom, V.I. Arnold showed the existence of a class of motions characterized by an unbounded ``secular change'' in the action variables, regardless of the perturbation size. Technically speaking, the proof relies on a \emph{mechanism}\footnote{the term ``Arnold's mechanism'' is often used in place of ``Arnold's diffusion''.} apt to exploit the absence of topological entrapment of invariant structures due to the problem dimension and connect neighbourhoods of invariant tori whose actions are ``slightly'' different one another. A diffusing set of trajectories follows as a consequence of the possibility to iterate this procedure over a suitably constructed sequence of invariant tori, called \emph{transition chain}. \\
The potential of this approach has been recognised and profitably developed in the extensive work \cite{AIHPA_1994__60_1_1_0}, motivated by the D'Alembert problem in Celestial Mechanics and dealing with a wide class of systems, including Arnold's. It is important to stress, for our purposes, that the proof of existence of diffusing trajectories of \cite{AIHPA_1994__60_1_1_0} possesses the form of an ``iterative algorithm''. It was immediately clear that those unstable trajectories were explicitly constructible.\\
However, the instability time $T_d$ estimates, appearing for the very first time in \cite{AIHPA_1994__60_1_1_0}, showed that such a time would have been extremely large: so large that the existence of those trajectories could have been extremely difficult (or even impossible) to be witnessed in concrete models (for a bound, compare \cite[$\S 8$]{AIHPA_1994__60_1_1_0} with \cite{AIHPA_1998__68_1_135_0}).\\    
Few years later, U. Bessi proposed a radically different approach with the paper  \cite{MR1375654}: whilst the techniques used in the proof contained  in  \cite{AIHPA_1994__60_1_1_0} (as those used in \cite{arn}) were mainly of geometric nature, Bessi's point of view was genuinely analytic and based on  Calculus of Variations tools. These results are to be compared with the well-known Mather's Theory, e.g. \cite{math}, still variational in nature but with several technical differences. See e.g. \cite{Xia1998} and \cite{bern} for developments of this Theory and \cite{bes2001} for more information on a possible comparison between the two approaches. \\
Coming back to the earliest Bessi's analysis, it was focused on the model 
\beq{eq:lag}
L_{\mu,\epsilon}(q,Q,\dot{q},\dot{Q},t)=\frac{\dot{Q}^2}{2}+\frac{\dot{q}^2}{2}+\epsilon (1-\cos q)-\mu \epsilon (1-\cos q)(\cos Q+\cos t) \mx{,}
\eeq
i.e. fully equivalent to the one considered \cite{arn}.\\
The key aspect was that, despite a step back with respect to the generality achieved in  \cite{AIHPA_1994__60_1_1_0} was taken, the bound on $T_d$ found by Bessi at this very first attempt was striking: ``just'' $O(\mu^{-2})$. The developments of this  ``variational school'', which flourished shortly thereafter, see e.g. \cite{MR1912262}, \cite{MR1996776} and other works, will improve the arguments in a way to obtain solutions for which $T_d=O(\mu^{-1} \log \mu^{-1})$. Knowingly, this bound has been shown to be optimal in \cite{MR1996776} for a class of systems including (\ref{eq:lag}). \\ 
So it was clear, since the first its first appearance in \cite{MR1375654}, that this approach was extremely powerful to prove the existence of ``fast'' drifting solutions. Despite ``why?'' could have been the most natural, but also the most difficult question to be answered, another question arose at that time, which certainly would have shed light on the first one: is there a way to use the method to compute those trajectories? (or, better, approximate them with a given precision), in other terms, is the variational method constructive? This problem, pointed out by G. Gallavotti, is the motivation of the present work. Clearly, this paper has not the conceit to be an answer to it, but it aims to be a step towards its understanding and it is not unlikely that the argument proposed could be easily extended to more general systems than the one studied here.

\section{Bessi's approach and the problem of constructivity}\label{sec:bessisapproach}
The idea behind Bessi's argument is recalled below for the reader's convenience and to fix the notations, by going along the lines of the presentation contained in \cite{MR1395239}. From now on, we shall refer to (\ref{eq:lag}) as the same system in which $\epsilon=1$ has been set, i.e. $L_{\mu}:=L_{\mu,1}$. We also recall that, in this case, the drift involves the variable $\dot{Q}(t)$, \cite{MR1375654}. Physically, it corresponds to an energy exchange (operated by the pendulum as a result of a multitude of small transfers) between the ``rotator'' (term $\dot{Q}^2/2$) and the harmonic oscillator (with frequency one), ``hidden'' in the time-dependent forcing. Clearly, if $\mu=0$ such a transfer would be impossible.\\
\smallskip \\
Let us start by noticing that, given $T_b>T_a$ and $Q_b>Q_a$, the boundary value problem  
\beq{eq:bvpuno}
\begin{cases}
\min \int_{T_a}^{T_b} L_{\mu}(q,\dot{q},Q,\dot{Q},t)dt \\
q(T_a)=(2 l-1)\pi \\  
q(T_b)=(2 l+1)\pi \\
Q(T_a)=Q_a\\
Q(T_b)=Q_b
\end{cases}
\eeq
has (in particular) a unique smooth solution for any $l \in \NN$ and any sufficiently small $\mu \geq 0$. The proof of this fact relies on standard tools of the Calculus of Variations, see e.g. \cite[P. 313]{MR1395239}. However, this information is not sufficient for our purposes and a quantitative criterion for the existence of such a minimum in a suitable discrete space will be given in Sec. \ref{sec:bvp}.  \\
Let us now fix $T_1<T_2<T_3$ and $Q_1<Q_2<Q_3$, then consider the two boundary value problems obtained from (\ref{eq:bvpuno}) by setting 
\[
(T_a,Q_a,T_b,Q_b,l)=(T_1,Q_1,T_2,Q_2,1), \qquad (T_a,Q_a,T_b,Q_b,l)=(T_2,Q_2,T_3,Q_3,2) \mx{,}
\]
and denoting with $(q^{(1,2)}(t),Q^{(1,2)}(t))$ and $(q^{(2,3)}(t),Q^{(2,3)}(t))$ their respective minimizers. Let us finally ``join'' the trajectories
\beq{eq:joincont}
(q^{(1,3)}(t),Q^{(1,3)}(t)):=
\begin{cases}
(q^{(1,2)}(t),Q^{(1,2)}(t)) & \qquad t \in [T_1,T_2]\\
(q^{(2,3)}(t),Q^{(2,3)}(t)) & \qquad t \in (T_2,T_3]\\
\end{cases}
\mx{.}
\eeq
If $\mu=0$ is easy to realise that $(q^{(1,3)}(t),Q^{(1,3)}(t))$ is never a minimiser for $\mathcal{F}:=\mathcal{F}_{1,2}+\mathcal{F}_{2,3}$ (see (\ref{eq:bvps}) below for their definition) unless $T_2 = (T_1+T_3)/2$. In fact, in the latter case, the two trajectories are indistinguishable from the unperturbed solution between $T_{1}$ and $T_{3}$ for $L_0$. Otherwise, if  $T_2 \neq (T_1+T_3)/2$, one has $\dot{q}^{(1,2)}(T_2) \neq \dot{q}^{(2,3)}(T_2)$ i.e. the trajectories speeds ``break'' at the joining point, and the same happens to $\dot{Q}$.\\
The described property has a well known interpretation from the splitting point of view: in the unperturbed system the only possible intersection is the homoclinic one, i.e., defined 
\beq{eq:omegai}
\omega_{i,i+1}:=\frac{Q_{i+1}-Q_i}{T_{i+1}-T_i} \mx{,}
\eeq
one can only have $\omega_{1,2}=\omega_{2,3}$ and such an intersection is \emph{flat}, see, e.g. \cite{for17}.\\
The situation is conceptually different if $\mu>0$, no matter how small, and this is the key ingredient of Bessi's argument. In fact, as it is shown in \cite{MR1375654}, it is possible to find $\{(T_i^*,Q_i^*)\}_{i=1,2,3}$ in such a way the obtained $(q^{(1,3)}(t),Q^{(1,3)}(t))$ \textbf{is} a minimiser for $\ml{F}$.\\
Knowingly, in order for this to happen, the values $\omega_{1,2}$ and $\omega_{2,3}$ have to be ``very close'', namely $\dist(\omega_{1,2},\omega_{2,3})=\ml{O}(\mu)$. This is strictly related, but not equivalent, see \cite{MR1996776}, to the closeness of the action of two consecutive tori in a transition chain, as originally used in \cite{arn}, and, once more, to their splitting properties.
\\ The key aspect is that the above described approach can be extended to the case of $N$ transitions. More precisely, by considering  $T_1<T_2<\ldots<T_{N+1}$ and $Q_1<Q_2<\ldots<Q_{N+1}$ and the family of BVPs
\beq{eq:bvps}
\min_{q^{(i,i+1)},Q^{(i,i+1)}} \mathcal{F}_{i,i+1} \mx{,}
\eeq
where
\[
\mathcal{F}^{(i,i+1)}=\int_{T_i}^{T_{i+1}} L_{\mu}(q^{(i,i+1)},Q^{(i,i+1)},\dot{q}^{(i,i+1)},\dot{Q}^{(i,i+1)},t)dt \mx{,}
\]
subject to 
\beq{eq:bv}
\begin{cases}
q^{(i,i+1)}(T_i)&=(2i-3)\pi \\ 
q^{(i,i+1)}(T_{i+1})&=(2i-1)\pi
\end{cases}
\qquad ; \qquad
\begin{cases}
Q^{(i,i+1)}(T_i)&=Q_i\\ 
Q^{(i,i+1)}(T_{i+1})&=Q_{i+1}
\end{cases} 
\mx{.}
\eeq
The task consists in finding a sequence $\{(T_i^*,Q_i^*)\}$ which minimises $\ml{F}:=\sum_{i=1}^N \ml{F}_{i,i+1}$.  Those point are shown to be ``close'' to a ``pseudo-sequence'' which plays the role of ``skeleton'' of the trajectory, constructed in such a way $\dist(\omega_i,\omega_{i+1})=\ml{O}(\mu)$ but $\dist(\omega_{1,2},\omega_{N,N+1})=\ml{O}(1)$. This is clearly achievable by choosing $N=\ml{O}(1/\mu)$, i.e. ``very large''. Furthermore, as there are no restrictions on $N$ (at least for the model at hand), such a sequence can be chosen arbitrarily long. In this way, by Lagrange's Theorem, there exist $\tilde{t}_1 \in (T_1,T_2)$ and $\tilde{t}_N \in (T_N,T_{N+1})$ such that $\dot{Q}^{(1,2)}(\tilde{t}_1)=\omega_{1,2}$ and $\dot{Q}^{(N,N+1)}(\tilde{t}_N)=\omega_{N,N+1}$ which is the desired result of drift.\\
It is possible to notice that the above described construction provides a sort of ``fictitious discretisation'' of the problem, in which the the action functional $\ml{F}$ can be interpreted as a function of $2(N-1)$ real variables:  
\beq{eq:hetfun}
\ml{F}(T_2,Q_2,\ldots,T_N,Q_N)=\sum_{i=1}^N \ml{F}_{i,i+1}(T_i,Q_i,T_{i+1},Q_{i+1}) \mx{,}
\eeq
Clearly, the previous sum controls on a ``substrate'' of functions $\ml{F}_{i+1}=\ml{F}_{i+1}(q^{(i,i+1)}(t),Q^{(i,i+1)}(t))$ in which the infinite dimensional nature of the problem is hidden. However, a large part of the theory can be carried out regardless of how the BVPs (\ref{eq:bvps}) are studied. This distinction into two \emph{layers} of problems, already used since the earliest work of Bessi, has remarkable consequences either for the proof of constructivity or for the efficiency of its implementation on a machine.\\
As a finite dimensional reduction of the classical action functional, we shall refer here to $\ml{F}$ as \emph{action function}. This object, already defined in \cite[Pag. 1129]{MR1375654}, has been used in a modified version (i.e. featuring conditions on the initial and final values of the actions) under the name of \emph{heteroclinic function} in \cite{MR1912262} and subsequent papers, see also \cite[Lemma 2.3]{MR1996776}. 
 \bigskip \\
Technically speaking, the problem of the constructivity of the variational approach, is motivated by the fact that, since \cite{MR1375654}, the minimum of the function $\ml{F}$ is shown to be attained in a compact set (more precisely the union of $N-1$ compact sets indexed by $i$, each of them containing $(T_i,Q_i)$) simply by proving that it cannot be on the boundary of such a set. Its existence follows directly by virtue of the Theorem of Weierstra\ss. Clearly, this does not guarantee the uniqueness of such a minimum as well as the possibility to approximate it via an algorithm.\\ 
As anticipated, a substantial difference with respect to the geometrical methods such as those used in \cite{AIHPA_1994__60_1_1_0} or \cite{Marco1996} (the latter also known as \emph{windows method}), is that the variational approach does not possess an intrinsic iterative structure and its constructivity cannot be established through a ``one-step-based'' proof.\\    
As for the last remark, it should be added that there are examples in the literature on geometrical methods, in which the iterative structure has been weakened or partly removed in order to achieve faster solutions and ``try to reach'' the diffusion speed obtained via variational tools. 
For instance, the refinement of the above mentioned windows method presented in \cite{cresguil3} claims to be ``hybrid'' in a certain sense and to have borrowed some ``ideas'' from the variational approach. Similarly, the argument used in \cite{for17} to overcome the correct alignment difficulties of the windows in the isochronous case, has only a partially iterative structure. Roughly speaking, the information used by iterative proofs is limited to the previous element of the transition chain, whilst the variational method, ``links'' the information of the whole trajectory. It is not surprising that some features, such as the speed of diffusion, are, in some sense,``optimised''. See Sec. \ref{sec:experiments} for further comments on this aspect.    
\medskip
\\
Another comment is in order concerning the concept of constructivity. Without addressing the deep and interesting concepts of the \emph{constructivism} and its implications, see e.g. \cite{sepmathconstructive}, one can say that a problem is defined as ``constructive'' if there exists a ``recipe'' \cite{goldstein2012constructive} to find (or, better, \emph{construct}) its solution. For instance, the classical Implicit Function Theorem can be classified as ``constructive'', as it can be proved by using an algorithm which is shown to be convergent to a (unique) solution by virtue of the Banach-Caccioppoli contraction principle. It should be added that this algorithm is generally ``very efficient'' and this feature will be profitably used in this paper. The proper ``constructivity'' of several results in Classical Analysis has often been argued, however, the proof by \emph{reductio ad absurdum} such as the one used in \cite{MR1375654} and \cite{MR1912262} stands amongst the almost unanimously recognised members of the ``non-constructive'' family, justifying in this way the problem discussed here. 
\medskip \\
It is important to stress that machine aided tools have a well established tradition in the Arnold's diffusion context. Computer based approaches have been profitably employed for the so called ``numerical evidence'' of diffusion, see e.g. \cite{MR2004513} and subsequent papers. 
For instance, in \cite{efthy}, an hybrid between analytic (normal forms) and numerical tools is used.\\
On the other hand, ``purely theoretical'' approaches may result in being either unsuitable or unsatisfactory in concrete applications. For this reason, numerical or computer assisted tools (or a mixture of both) have successfully took part in the context of rigorous proofs, see e.g. \cite{MR3551192}. Another case of computer assisted proof is contained in \cite{CAPINSKI2022105970}. \\
Clearly, the mentioned works are just a taster of the wide literature on this subject and even an essential review would go far beyond the purposes of this paper. On the other hand, the spirit of this work is quite different as the diffusion result is known, but not any detail, to the best of our knowledge, can be found (to date) in the literature about its constructivity. It has to be stressed that the problem is not limited to the possibility of an implementation of the method but involves, at the same time, an estimation of the computational cost of the whole task.

\section{Discretisation of the action functional and main result}\label{sec:discr}
\begin{figure}[htbp!]
\centering
\includegraphics[scale=0.85]{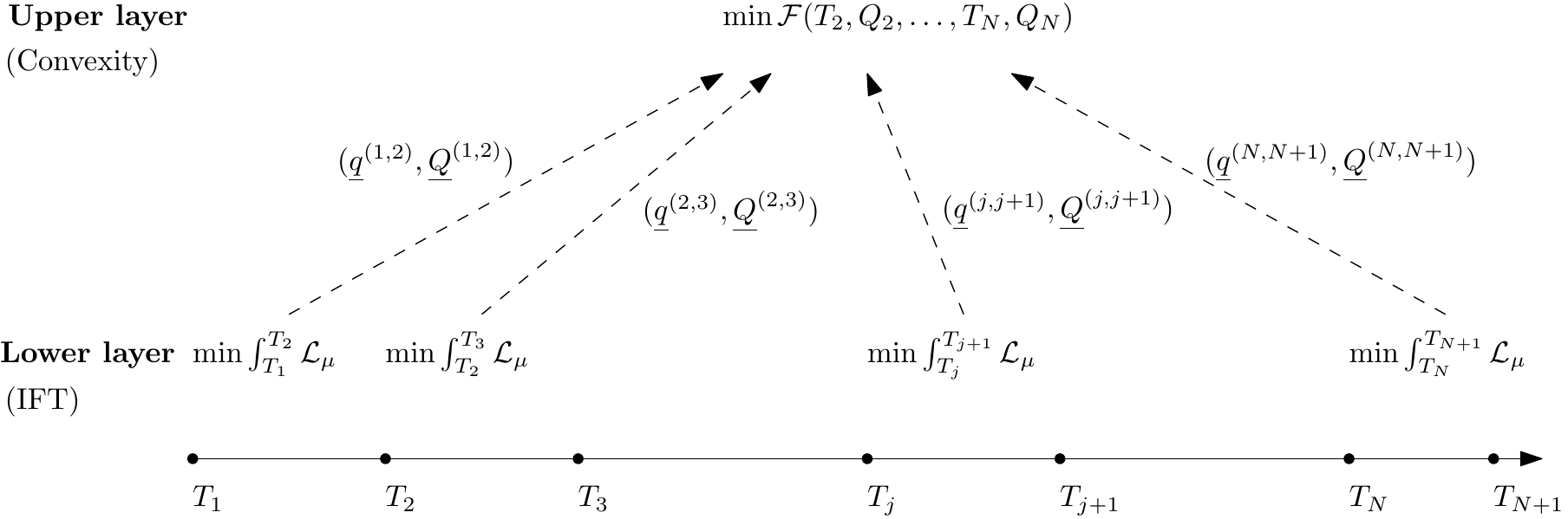}
\caption{The multilayer structure of the proof / implementation.}
\label{fig:layers}
\end{figure}

As we are interested in constructing a concrete approximation of the diffusing trajectories we have to deal with a standard finite dimensional reduction obtained by discretisation. For (much) more details about this approach, its implications and a broad panorama of the field of Variational Integrators we refer to the comprehensive work \cite{mawe}.\\
Let us suppose for a moment that the number of transitions $N$ and a ``sufficiently large'' $n \in \NN$  have been fixed. Then we consider, for all $i=1,\ldots,N$, the set of ``internal nodes'' of $[T_i,T_{i+1}]$
\beq{eq:part}
M_{i,i+1}:=\{t_k^{(i)}=T_i+h_i k,\, k=1,\ldots,n_i \}, \qquad (n_i+1) h_i:=T_{i+1}-T_i \mx{,}
\eeq
(in particular $t_0^{(i)}\equiv T_i$, $t_{n_i+1}^{(i)}\equiv T_{i+1}$), supposing $\min_{i=1,\ldots,N}\{n_i\} \geq n$. Let us preliminarily observe that the notational setting introduced below can be straightforwardly adapted \emph{mutatis mutandis} if $M_{i,i+1}$ is replaced with $M_{i,i+1}^*:=\{T_i\} \cup M_{i,i+1} \cup \{T_{i+1} \}$. It is easy to realize, however, that the latter has a marginal relevance as the values of the solutions at the extrema of $[T_i,T_{i+1}]$ are known from the boundary conditions. \medskip \\
Let us start by defining $\mathfrak{D}(M_{i,i+1})$ as the set of all the \emph{discrete functions} with domain $M_{i,i+1}$. Hence, an element of this set is a $n_i-$ple  $\ue{f} \equiv (f_1,\ldots,f_{n_i})$, where $f_j$ is interpreted as the value of $\ue{f}$ at the time $t_j^{(i)}$. Clearly we can identify $M_{i,i+1} \ni \ue{t}^{(i)}:=(t_1^{(i)},\ldots,t_{n_i}^{(i)})$. On the other hand, given a function $g:[T_i,T_{i+1}] \rw \RR$, it is straightforward to associate this notation to the \emph{restriction} of $g$ to $M_{i,i+1}$ and we shall set $\ue{g}:=g|_{M_{i,i+1}}$ where $g_k:=g(t_k^{(i)})$ for all $k=1,\ldots,n_i$.\\ 
For all $i$, the set $M_{i,i+1}$ will be endowed with the standard infinity norm $\norm{\cdot}$ i.e. $\norm{\ue{f}}:=\max_{k=1,\ldots,n_i}|f_k|$, for all $\ue{f} \in \mathfrak{D}(M_{i,i+1})$. It is immediate to verify that $(\mathfrak{D}(M_{i,i+1}),\norm{\cdot})$ is a Banach space.\\
Given $\ue{f}\in \mathfrak{D}(M_{i,i+1})$ and $\ue{g}\in \mathfrak{D}(M_{i+1,i+2})$, we define 
\beq{eq:union}
\ue{f} \cup \ue{g} := (f_1,\ldots,f_{n_i},g_{1},\ldots,g_{n_{i+1}}) \in \mathfrak{D}(M_{i,i+1}) \cup \mathfrak{D}(M_{i+1,i+2}) \mx{.}
\eeq
It is immediate to extend inductively the previous definition to a finite number of elements. 
\medskip \\
Let us now deal with the discretisation of $\ml{F}$. Following a widely used approach, we write $\ml{F}_{i,i+1}=\sum_{k=1}^{n_i} \int_{t_k}^{t_{k+1}} L_{\mu}$ then we compute each integral via the \emph{rectangle rule}. In this way 
\beq{eq:discretef}
\ml{F}_{i,i+1} \sim \ml{F}_{i,i+1}|_{M_{i,i+1}} := h_i \sum_{k=0}^{n_i} L_k^{(i,i+1)} \mx{,}
\eeq
where
\beq{eq:discretel}
L_k^{(i,i+1)}:=L(q_k^{(i,i+1)},h_i^{-1}(q_{k+1}^{(i,i+1)}-q_k^{(i,i+1)}),Q_k^{(i,i+1)},h_i^{-1}(Q_{k+1}^{(i,i+1)}-Q_k^{(i,i+1)}),t_k) \mx{,} 
\eeq
and the values $q_{0}^{(i,i+1)},q_{n_i+1}^{(i,i+1)},Q_{0}^{(i,i+1)}$ and $Q_{n_i+1}^{(i,i+1)}$ are given by (\ref{eq:bv}).\\
The use of the rectangle rule, together with the ``simple'' forward finite difference formula for the derivative 
\beq{eq:discrete}
\dot{q}_k^{(i,i+1)}\sim h_i^{-1}(q_{k+1}^{(i,i+1)}-q_k^{(i,i+1)}), \qquad 
\dot{Q}_k^{(i,i+1)}\sim h_i^{-1}(Q_{k+1}^{(i,i+1)}-Q_k^{(i,i+1)})\mx{,} 
\eeq
appearing in (\ref{eq:discretel}), can be certainly accounted amongst the simplest possible approaches in this context. \\
The quantity $L_k^{(i,i+1)}$, often referred to as \emph{discrete Lagrangian}, see e.g. \cite{mawe}, depends on two unknowns $\ue{q}^{(i,i+1)},\ue{Q}^{(i,i+1)}\in \mathfrak{D}(M_{i,i+1})$ which have the meaning of ``approximated trajectory'' for $\ml{F}_{i.i+1}$.\\We are now interested in studying the following problem 
\beq{eq:minproblem}
\min_{\{T_i,Q_i\}_{i=2,\ldots,N}} \ml{F}|_{M_{1,N+1}}, \qquad \ml{F}|_{M_{1,N+1}}:=\sum_{i=1}^{N} \min_{\left(\ue{q}^{(i,i+1)},\ue{Q}^{(i,i+1)}\right)} \ml{F}_{i,i+1}|_{M_{i,i+1}} \mx{,}
\eeq
where $M_{1,N+1}:=\cup_{i=1}^N M_{i,i+1}$ and $T_1,Q_1,T_{N+1},Q_{N+1}$ are interpreted as fixed, as anticipated in Sec. \ref{sec:bessisapproach}. It is natural to define $\mathfrak{D}(M_{1,N+1}):=\cup_{i=1}^N \mathfrak{D}(M_{i,i+1})$. \\
Hence, the solution of (\ref{eq:minproblem}) will be an object of the form  
\beq{eq:joindiscr}
(\ue{q}^{(1,N+1)},\ue{Q}^{(1,N+1)})=\bigcup_{i=1}^N (\ue{q}^{(i,i+1)},\ue{Q}^{(i,i+1)}) \in \mathfrak{D}(M_{1,N+1}) \mx{,}
\eeq
(recall (\ref{eq:union})), i.e. ``joining'' all the ``pieces of trajectory'' indexed by $i$. Either $\ue{q}^{(1,N+1)}$ or $\ue{Q}^{(1,N+1)}$ have $\ml{N}:=\sum_{i=1}^N n_i$ ``components''. \\
Let us finally introduce the following notation: given a subset $\mathcal{U} \subset \RR^n$, we define $\diam_-(\mathcal{U}):=\sup\{d>0:B_{d/2}(x_0)\subseteq \mathcal{U}, \, \forall x_0 \in \mathcal{U}\}$.\\
In this setting, the main result states as follows
\begin{satz}\label{thm:main} 
Suppose $\mu$ ``sufficiently small'' (more on this later) and define  
\beq{eq:tmenopiu}
\ml{T}^-(\mu):=(3/4)\log(320/\mu),  \qquad \ml{T}^+(\mu):=\pi\log(640/\mu)  \mx{.}
\eeq
Then there exists (at least) an interval $\chi_{d,\mu} \subset [0,1]$ with $\diam_-(\chi_{d,\mu})=O(1)$ such that, for any $\omega_I, \omega_F \in \chi_{d,\mu}$ with $\omega_I < \omega_F$, it is possible to construct a sequence $\left\{ \left(  \ue{q}^{r,(1,N+1)}, \ue{Q}^{r,(1,N+1)} \right)_r \right\} \in \mathfrak{D}(M_{1,N+1})$ (i.e.  of objects of the form (\ref{eq:joindiscr})), with 
\beq{eq:choicen}
N=4+2\lceil (\omega_F-\omega_I) /(2 C \mu)\rceil \mx{,} 
\eeq
and $C:=1/20$, 
which converges to the unique absolute minimum of $\ml{F}|_{M_{1,N+1}}$ on $\mathfrak{D}(M_{1,N+1})$ for $r \rw \infty$,  with the following properties: 
\begin{enumerate}
\item There exist $\delta=O(1)$ and $(\tilde{n}_1,\tilde{n}_N) \in [2,n_1-1]\times [2,n_N-1] \subset \NN^2$ such that 
\beq{eq:shadow}
\max\{|\dot{Q}_{\tilde{n}_1}-\omega_I |,|\dot{Q}_{\tilde{n}_N}-\omega_F |\} \leq  (4 \delta /\ml{T}^-(\mu)+\pi \mu)+o(h,\mu)\mx{,} 
\eeq
where 
$h:=\max_i h_i$, see (\ref{eq:part}). Note that the time derivatives appearing in (\ref{eq:shadow}) are understood in the sense of (\ref{eq:discrete}). 
\item The time of ``drift'', i.e. $\mathcal{T}_d:=T_{N+1}-T_1$, is bounded as follows
\beq{eq:drifttime}
\mathcal{T}_d \leq N \mathcal{T}^+(\mu)\mx{.}
\eeq
\end{enumerate}
\end{satz}
First of all, we notice that the time estimate obtained agrees with the well known $T_d =\ml{O} (\mu^{-1}\log \mu^{-1})$, which we recall to be optimal, as shown in \cite{MR1996776}.\\ 
The bound (\ref{eq:shadow}) represents the arbitrary separation between the initial and final value of $\dot{Q}$, which is the main implication of the Arnold's diffusion phenomenon. The constant $\delta$ appearing in (\ref{eq:shadow}) is defined in Sec. \ref{sec:convexity}, and represents the size of the ``boxes'' where the $(T_i,Q_i)$ are located. However, as it will be clear, $\delta$ can be reduced at a price of increasing the overall time of transition, in other words, a ``more precise'' shadowing requires ``more'' time. \\
Clearly, the bound (\ref{eq:shadow}) is not particularly expressive if $\omega_I$ and $\omega_F$ are close to each other, more precisely, every time $\omega_F-\omega_I=O(\mu)$. However, this should not be seen as a limitation of such a bound but as a general feature of the solutions obtained with this method. In fact, despite the argument forces $\dot{Q}$ to increase over a $O(1/\mu)$ time-span, there is no restriction for such a value to ``decrease'' or ``oscillate'' over $O(1)$ time intervals. This is most likely the ``secret'' behind the optimality (in time) of the trajectories constructed in this way. More on this aspect in Sec. \ref{sec:experiments}. \\
\medskip \\
The setting of Sec. \ref{sec:bessisapproach} has made clear that the required machine implementation is nothing but a large-scale optimisation problem, in fact, it is immediate to realise that the dimension of the discrete space involved is roughly $\ml{N}=\ml{O}(n N)$, where $n$ is typically ``very large''. As we shall see in Sec. \ref{sec:experiments}, fairly interesting diffusion phenomena could easily involve $\ml{N}=O(10^{8})$ variables and this clearly requires a strategy to tackle the related CPU time. \\   
First of all, let us notice that the solution of each BVP (\ref{eq:bvps}) indexed by $i$, is fully determined once $T_i,T_{i+1},Q_i,Q_{i+1}$ have been set. If we suppose for a moment that we can determine $\mu_0$ such that $\ml{F}_{i,i+1}|_{M_{i,i+1}}$ possesses a (unique) absolute minimum for any (suitable) choice of $T_i,T_{i+1},Q_i,Q_{i+1}$ and \underline{all} $\mu \in (0,\mu_0]$, we can treat $\ml{F}$ as a function of $(T_2,Q_2,\ldots,T_N,Q_N)$ ``only'', and the minimisation of $\ml{F}$ is reduced from a problem in dimension $\RR^{2\ml{N}}$ to one in dimension $\RR^{2(N-1)}$.\\
The described dimensional reduction exploits the nature of the two different layers of problems. In particular, as it will be clear later, the minimisation of $\ml{F}_{i,i+1}$ on $\RR^{2 n_i}$, \emph{lower layer}, possesses an ``implicit function'' structure. In this way, the classical Implicit Function Theorem (IFT) can be used either to obtain a rigorous proof of existence and uniqueness of the problem (\ref{eq:bvpuno}), or to approximate such a minimum on a machine very efficiently. \\
On the other hand, the minimisation of $\ml{F}$ on $\RR^{2(N-1)}$, \emph{upper layer}, does not possess such a strong structure and only ``weaker'' arguments, such as convexity, which can be proven up to a first order in $\mu$ (hence the ``sufficiently small $\mu$'' in the statement). Fig. \ref{fig:layers} gives a pictorial idea of how the whole minimization process is organised. \\
\\Given the algorithmic structure of the IFT, the existence of any valid algorithm for the upper layer minimisation solves the problem of the constructivity. The class of methods commonly known as ``Nesterov accelerated gradient'' (appeared, in its first ``variant'', in \cite{Nesterov1983AMF}), turns out to be particularly adequate in our case. Firstly, it is a gradient-based method, and in our case  
the evaluation of the gradient is exact, as, remarkably, its explicit expression is known:
\begin{align}
\pl_{T_i} \ml{F}& = 2^{-1}\left(\dot{q}_{i,i+1}^2(T_i)-\dot{q}_{i-1,i}^2(T_i) \right) + 2^{-1} \left( \dot{Q}_{i,i+1}^2(T_i)-\dot{Q}_{i-1,i}^2(T_i)\right) \mx{,} \label{eq:gradbertiuno}\\
\pl_{Q_i} \ml{F} & =\dot{Q}_{i-1,i}(T_i)- \dot{Q}_{i,i+1}(T_i) \label{eq:gradbertidue} \mx{,}
\end{align}
see \cite[Lemma 2.2]{MR1996776} for the proof. The latter explicit expression makes it cost effective as well, as the boundary values are known from the (approximate) solutions on the lower layer.  \\
Moreover, the convergence of such a method follows from known results (more on this aspect later), once the (strict) convexity of $\ml{F}$ has been shown. Knowingly, the rate of convergence of Nesterov's method is $O(r^{-2})$, where $r$ is the step-number, differently from the classical gradient descent, in which the rate is estimated as $O(r^{-1})$. This offers a remarkable advantage in our case as the evaluation of a single step is generally ``very costly''.\\
In any case, it has to be stressed that our aim is to show the full constructivity of the setting, and this is done by providing the simplest approach in terms of implementation. The improvements that could be introduced to our implementation are clearly countless, but they are just beyond the purposes of this paper.
\smallskip \\
As anticipated in the abstract, the proof contained in this paper relies on already existing ideas, in particular on those contained in the papers \cite{MR1375654} and \cite{MR1996776}, but some additional ingredients are necessary in order to conclude on the constructivity problem. For instance, the statements in Secs. \ref{sec:three} and \ref{sec:convexity} are formulated in terms of Jacobian elliptic functions for the pendulum. In this way, it is possible to use the well established theory behind them in order to obtain key quantitative estimates. At the same time, from the implementation point of view, elliptic functions such as $\am(u,\kappa)$, which naturally appear in the solution of the unperturbed pendulum, can be efficiently computed via well known approximation formulae or, better, via software in-built functions.\\
The formulation in terms of Jacobian elliptic functions leads to an object which is immediately recognized in the first order approximation of the action functional, and it can be naturally interpreted as an ``extension'' of the standard Melnikov integral outside the pendulum separatrix, see Prop \ref{prop:extended}. The latter, is shown to play a key role in the proof of the dominance of such a first order term for sufficiently small $\mu$. \\
Other tools are needed in order to deal with the discrete setting of the result. More precisely, the already mentioned use of the classical IFT for the BVP (\ref{eq:bvpuno}), requires a non-trivial effort in order to prove the non-degeneracy of the Jacobian of the system. \\
\smallskip \\
The paper is organized as follows: In Sec. \ref{sec:three} the first order expansion of the action function $\ml{F}$ in terms of elliptic functions and the mentioned ``extended Melnikov'' integral are derived. In Sec. \ref{sec:convexity} we show the possibility to approximate the minimum of $\ml{F}$ (upper-layer constructivity) by showing the strict convexity of $\ml{F}$ under standard assumptions. This is achieved by assuming the constructivity on the lower-layer, which will be thereafter shown in Sec. \ref{sec:bvp}. \\
The final section is devoted to the discussion of the algorithm implementation and the aspects related to its performance, alongside the outcomes of some of the most relevant simulations.

\section{Approximation of the action function}\label{sec:three}
In this section we obtain some approximation formulae that will be used afterwards to obtain the required first order approximation of the action function. The essential notions from the theory of Jacobian elliptic functions used here are briefly recalled in Appendix A for the reader's convenience. Extensive expositions can be found, e.g., in \cite{abramowitz1965handbook} and \cite{whittaker1996course}.  
\begin{prop} The following expansion holds
\[
\ml{F}=\ml{F}^{(0)}+\mu \ml{F}^{(1)}+o(\mu) \mx{,}
\]
where 
\begin{align}
\ml{F}^{(0)}&= 2 \sum_{i=1}^N \left\{ k_{i,i+1}^{-1} \left[2 \E{k_{i,i+1}} +(k_{i,i+1}')^2 \K{k_{i,i+1}}\right]+\omega_{i,i+1}^2 T_{i,i+1}^- \right\} \label{eq:fzero}\\
\ml{F}^{(1)}&= - \Gamma_{1,2}^+(\omega_{1,2}) \cos Q_1 - \Gamma_{1,2}^+(1) \cos T_1+\Theta_{1,2}^+(\omega_{1,2}) \sin Q_1 +\Theta_{1,2}^+(1) \sin T_1 \nonumber\\
& -\sum_{i=2}^N \biggl\{ \left[\Gamma_{i-1,i}^-(\omega_{i-1,i})+\Gamma_{i,i+1}^+(\omega_{i,i+1})\right] \cos Q_i  + \left[\Gamma_{i-1,i}^-(1)+\Gamma_{i,i+1}^+(1)\right] \cos T_i \nonumber \\
& -\left[\Theta_{i-1,i}^-(\omega_{i-1,i})+\Theta_{i,i+1}^+(\omega_{i,i+1})\right] \sin Q_i  -\left[\Theta_{i-1,i}^-(1)+\Theta_{i,i+1}^+(1)\right] \sin T_i \biggl\} \nonumber \\
& - \Gamma_{N,N+1}^-(\omega_{N,N+1}) \cos Q_{N+1} - \Gamma_{N,N+1}^-(1) \cos T_{N+1} \nonumber \\
&+\Theta_{N,N+1}^-(\omega_{N,N+1}) \sin Q_{N+1} +\Theta_{N,N+1}^-(1) \sin T_{N+1} \label{eq:funo}
\end{align}
with
\begin{align}
\Gamma_{i-1,i}^-(\omega)&:=2\int_{-T_{i-1,i}^-}^0 \cn^2(k_{i-1,i}^{-1}t,k_{i-1,i})\cos(\omega t) dt, \label{eq:gammam}\\ 
\Theta_{i-1,i}^-(\omega)&:=2\int_{-T_{i-1,i}^-}^0 \cn^2(k_{i-1,i}^{-1}t,k_{i-1,i})\sin(\omega t) dt, \label{eq:thetam} \\
\Gamma_{i,i+1}^+(\omega)&:=2\int_{0}^{T_{i,i+1}^- } \cn^2(k_{i,i+1}^{-1}t,k_{i,i+1})\cos(\omega t) dt, \label{eq:gammap}\\ 
\Theta_{i,i+1}^+(\omega)&:=2\int_{0}^{T_{i,i+1}^-} \cn^2(k_{i,i+1}^{-1}t,k_{i,i+1})\sin(\omega t) dt\label{eq:thetap} \mx{,}
\end{align}
and $T_{i-1,i}^{\pm}:=(T_i \pm T_{i-1})/2$ for all $i$. The variables $k_{i-1,i},k_{i,i+1}$ are computed in such a way, for all $i=2,\ldots,N$, see (\ref{eq:period}), one has 
\beq{eq:k}
T_{i-1,i}^-=k_{i-1,i} \K{k_{i-1,i}}, \qquad T_{i+1,i}^-= k_{i,i+1} \K{k_{i,i+1}} \mx{.}
\eeq
We shall denote $\mathcal{F}^{(\leq 1)} :=\mathcal{F}^{(0)}+\mu \mathcal{F}^{(1)}$.
\end{prop} 
\proof
First of all we notice that the following approximation holds
\beq{eq:firstorderapp}
\begin{array}{rcl}
\de \int_{T_a}^{T_b} L_{\mu}(q,\dot{q},Q,\dot{Q},t)dt & = & \de \int_{T_a}^{T_b} L_{0}(q_0,\dot{q}_0,Q_0,\dot{Q}_0,t)dt\\
&-& \de \mu \int_{T_a}^{T_b}
(1-\cos q_0(t))(\cos Q_0(t)+\cos t)dt 
 +o(\mu) \mx{.}
\end{array}
\eeq
This is easily shown by writing $q(t)=q_0(t)+v(t;\mu)$ and $Q(t)=Q_0(t)+w(t;\mu)$ where $v(T_a)=v(T_b)=w(T_a)=w(T_b)=0$ then expanding up to the first order in $\mu$. The required formula follows after an integration by parts (in which the boundary conditions for $v$ and $w$ are used) and from the minimality of $q_0(t)$ and $Q_0(t)$. \\
By using (\ref{eq:firstorderapp}) and the explicit form of $Q_0(t)$, one gets
\begin{align}
\ml{F}&=\sum_{i=1}^N \int_{T_i}^{T_{i+1}} L_0(q_0^{(i,i+1)}(t),\dot{q}_0^{(i,i+1)}(t),\dot{Q}_0^{(i,i+1)}(t))dt \nonumber \\ 
& -\mu \Biggl\{ \int_{T_1}^{T_{1,2}^+} (1-\cos q_0^{(1,2)}(t))(\cos (Q_{1}+\omega_{1,2}(t-T_{1}))+\cos t)dt \nonumber \\
&+\sum_{i=2}^N \left[ \int_{T_{i-1,i}^+}^{T_i} (1-\cos q_0^{(i-1,i)}(t))(\cos (Q_{i-1}+\omega_{i-1,i}(t-T_{i-1}))+\cos t)dt \right.\nonumber \\
& \left. +\int_{T_{i}}^{T_{i,i+1}^+} (1-\cos q_0^{(i,i+1)}(t))(\cos (Q_{i}+\omega_{i,i+1}(t-T_{i}))+\cos t)dt \right] \nonumber \\
& + \int_{T_{N,N+1}^+}^{T_N} (1-\cos q_0^{(N,N+1)}(t))(\cos (Q_{N}+\omega_{N,N+1}(t-T_{N}))+\cos t)dt \Biggl \} + o(\mu) \mx{.}\nonumber
\end{align}
Let us now focus our attention on $\ml{F}^{(0)}$ i.e. the first sum appearing in the previous expression. By performing suitable time translations we have 
\[
\ml{F}^{(0)}=2 \sum_{i=1}^N k_{i,i+1}^{-1} \int_{-\K{k_{i,i+1}}}^{\K{k_{i,i+1}}} \tilde{L}_{i,i+1}(\tau)d \tau \mx{,}
\]
where, see (\ref{eq:lagraell}),  
\[
\tilde{L}_{i,i+1}(\tau):=\dn^2(\tau,k_{i,i+1})+k_{i,i+1}^2 \cn^2 (\tau,k_{i,i+1})
\]
is proportional to the time rescaled Lagrangian of the unperturbed pendulum.\\
In this formulation, the integrals appearing in $\ml{F}^{(0)}$ are immediate from the theory of Jacobian elliptic functions, see e.g. \cite[Sec. 16.26]{abramowitz1965handbook}. Namely, by dropping the indices for simplicity,  
\[
\int_{-\K{k}}^{\K{k}} \tilde{L}(\tau) d \tau= \int_{-\K{k}}^{\K{k}} \dn^2(\tau,k)+k^2 \int_{-\K{k}}^{\K{k}} \cn^2(\tau,k)d \tau =4\E{k}+2(k')^2 \K{k} \mx{.}
\]
As for $\ml{F}^{(1)}$, by a similar set of time translations, one gets   
\begin{align}
\ml{F}^{(1)}&= \int_{0}^{T_{1,2}^-} (1-\cos \tilde{q}_0^{(1,2)}(\tau))
\left[\cos (Q_{1}+\omega_{1,2}\tau)+\cos (T_{1}+\tau) \right]d \tau \nonumber\\
&-\sum_{i=2}^N \biggl\{ \int_{-T_{i-1,i}^-}^{0} (1-\cos \tilde{q}_0^{(i-1,i)}(\tau)) \left[\cos (Q_{i}+\omega_{i-1,i}\tau)+\cos (T_i+\tau) \right]d \tau \nonumber \\
&  +\int_{0}^{T_{i,i+1}^-} (1-\cos \tilde{q}_0^{(i,i+1)}(\tau))
\left[\cos (Q_{i}+\omega_{i,i+1}\tau)+\cos (T_{i+1}+\tau) \right]d \tau
\biggl\} \nonumber\\
& \int_{-T_{N,N+1}^-}^{0} (1-\cos \tilde{q}_0^{(N,N+1)}(\tau)) \left[\cos (Q_{N+1}+\omega_{N,N+1}\tau)+\cos (T_{N+1}+\tau) \right]d \tau
\label{eq:funolast}
\end{align}
where 
\beq{eq:unp}
\tilde{q}_0^{(m-1,m)}(\tau)=\pi + 2 \am (k_{m-1,m}^{-1} \tau,k_{m-1,m})
\eeq
are the solutions of the unperturbed pendulum satisfying, for all $i=2,\ldots,N$,
\[
\tilde{q}_0^{(i-1,i)}(-T_{i-1,i}^-)=0,\qquad \tilde{q}_0^{(i-1,i)}(0)=\pi=\tilde{q}_0^{(i,i+1)}(0),\qquad \tilde{q}_0^{(i,i+1)}(T_{i,i+1}^-)=2 \pi \mx{.}
\] 
We recall that $\dot{\tilde{q}}_0^{(i-1,i)}(0) \neq \dot{\tilde{q}}_0^{(i,i+1)}(0)$ unless $T_{i-1,i}=T_{i,i+1}$ i.e. $T_i=(T_{i-1}+T_{i+1})/2$.\\
It is now sufficient to use (\ref{eq:unp}) in (\ref{eq:funolast}) and define the quantities (\ref{eq:gammam}), (\ref{eq:thetam}), (\ref{eq:gammap}) and (\ref{eq:thetap}), in order to get (\ref{eq:funo}).  
\endproof
\noindent The following statement will play a key role in the next section.
\begin{prop}\label{prop:extended} The following formula holds
\beq{eq:resint}
\int_{-k \K{k}}^{k \K{k}} \cn^2(k^{-1}t,k) \cos(\omega t) dt  =  \de \frac{\pi \omega}{ \sinh(\omega k \Kp{k})} + \ml{R}(k) \mx{,}
\eeq
where, recalling $\sd(u,k):=\sn(u,k)/\dn(u,k)$, 
\beq{eq:resto}
\ml{R}(k):= -(k')^2
\frac{4 k \Kp{k} \sin (\omega k \K{k})}{1-e^{-2 \omega k \Kp{k}}}
\int_0^1 \sd^2(2 s \Kp{k},k') e^{-2 \omega k s \Kp{k}} ds \mx{.}
\eeq
\end{prop}
\proof After a time rescaling $t=:k \tau$, (\ref{eq:resint}) is obtained by contour integration over the rectangle $ABCD$ of the complex plane, where $A=-\K{k}$, $B=\K{k}$, $C=\K{k}+2 i \Kp{k} $ and $D=-\K{k}+2i \Kp{k}$ then using 
\[
\Res \left[\cn^2(\tau,k)e^{i \omega k \tau},\tau=i \Kp{k} \right]=-k^{-1} i \omega e^{-\omega k \Kp{k}} \mx{.}
\]
The latter follows  from the expansion $\cn(z + i \Kp{k})=-(k z)^{-1}i+(6 k)^{-1}(2 k^2-1)i z +\ml{O}(z^2)$, see \cite[Pag. 504]{whittaker1996course}.  
\endproof

\begin{rem}
Due to the presence of $k'$ and otherwise bounded functions, the ``remainder'' $\ml{R}(k)$ vanishes as $k \rw 1$. Furthermore $\Kp{1}=\pi/2$ and $\cn (u,1)=\sech u$. Hence, by multiplying both sides of (\ref{eq:resint}) by $2$, its limit as $k \rw 1$, is nothing but the celebrated Melnikov integral 
\beq{eq:melnikov}
\ml{M}(\omega):=\int_{-\infty}^{+\infty} \sech^2 \tau \cos (\omega \tau)=\frac{ 2\pi \omega}{\sinh(\omega \pi /2)} \mx{.}
\eeq
Formula (\ref{eq:resint}) suggests that, by choosing $-T_a=T_b=k \K{k}$ then using the time $\tau$, the expansion of the action functional (\ref{eq:firstorderapp}), is explicitly computable up to the first order also \underline{outside} the separatrix, i.e. for $k <1$. This offers a great advantage in this context, as the so-called Melnikov approximation generally requires careful estimates of the difference between the exact action functional and the Melnikov term, see e.g. \cite{MR1375654}. 
\end{rem}

\section{Convexity and existence of the minimum}\label{sec:convexity}
Our aim is now to use the approximation formulae found in the previous section in order to show that the action function has a unique minimum on a suitable set, and this will be proven by convexity. We stress that the results of this section are independent of where the BVPs are studied, more precisely, they hold either in the infinite dimensional setting (\ref{eq:bvps}) or in their finite dimensional reduction. \\
Let us set for brevity $Z:=(T_2,Q_2,\ldots,T_N,Q_N)$. Given an arbitrary double sequence of indices $\ml{I}:=\{(\beta_1^{(i)},\beta_2^{(i)})\}_{i=2,\ldots,N}$, with $\beta_j^{(i)}\in \NN$, which is assumed to be strictly increasing, i.e. $\beta_j^{(i+1)}>\beta_j^{(i)}$, $j=1,2$, let us define
\[
\mathfrak{C}_{\ml{I},\delta}:=\{Z \in \RR^{2(N-1)} : (T_i,Q_i) \in [2 \pi \beta_1^{(i)} - \delta, 2 \pi \beta_1^{(i)} + \delta ]\times [2 \pi \beta_2^{(i)} - \delta, 2 \pi \beta_2^{(i)} + \delta ] \} \mx{.}
\]
\begin{rem}\label{rem:convexity}%
Let us consider the function $\mathfrak{m}:\RR \mapsto [-\pi,\pi]$ defined by $\mathfrak{m}(x):=(x+\pi)\Mod (2 \pi)-\pi$ and two increasing sequences $\{T_i^0\}$ and $\{Q_i^0\}$ satisfying $\mathfrak{m}(Z_i^0)\in [-\delta,\delta]$. It follows immediately from the definition that $Z \in \mathfrak{C}_{\ml{I},\delta}$ for some suitable $\ml{I}^0$. Any confusion should be avoided between the point 
$(T_i^0,Q_i^0)$ and the variables $(T_i,Q_i)$ which vary in the connected component of the set (i.e. the square of side $2 \delta$) containing $(T_i^0,Q_i^0)$.
\end{rem}
The main result of this section states as follows
\begin{lem}\label{lem:convexity} Let $C>0$ be a constant to be determined and let us assume the following three hypotheses
\begin{enumerate}
\item ``sufficiently large'' transition times
\beq{eq:hypuno}
k_{i,i+1} \geq k_0:=(1+C \mu)^{-\frac{1}{2}}, \qquad \forall i=1,\ldots,N.
\eeq
\item ``sufficiently close'' consecutive frequencies
\beq{eq:hypdue}
|\omega_{i,i+1}-\omega_{i-1,i}| \leq C \mu \mx{,}
\eeq
\item range of frequencies
\beq{eq:hyptre}
\omega_{i,i+1} \in (0,1], \qquad \forall i=1,\ldots,N.
\eeq
\end{enumerate}
Suppose in addition that, for any $T_{i+1}>T_i$ and $Q_{i+1}>Q_i$, it is possible to find $\mu_0$ such that, for any $\mu \in (0,\mu_0]$, there exist a unique absolute minimum $(\ue{q}^{(i,i+1)},\ue{Q}^{(i,i+1)}) \in \mathfrak{D}(\ml{M}_{i,i+1})$ of $\ml{F}_{i,i+1}|_{\ml{M}_{i,i+1}}$, for all $i=1,\ldots,N$. Let us finally set $\bar{\mu}:=10^{-5}$.\\
Then, for sufficiently small $\mu \leq \min\{\bar{\mu},\mu_0\}$, $C=1/20$ and all $\ml{I}$, the action function $\ml{F}$ possesses a unique absolute minimum  
\[
Z^* \in \mathfrak{C}_{\ml{I},\pi/4} \setminus \pl \mathfrak{C}_{\ml{I},\pi/4}\mx{,}
\] 
hence, we have set $\delta:=\pi/4$. Furthermore, for all $Z_0 \in \mathfrak{C}_{\ml{I},\pi/4} \setminus \pl \mathfrak{C}_{\ml{I},\pi/4}$ and all $\alpha \leq 1/2$, one has $Z^*=\lim_{r \rw \infty} Z_{r}$, where $Z_r$ is the ``Nesterov sequence'' defined by, see, e.g. \cite{JMLR:v17:15-084}
\beq{eq:graddesc}
\begin{aligned}
W_{r+1}&=Z_{r}-\alpha \nabla_{Z} \ml{F}(Z_r) \\
Z_{r+1}& = W_{r+1}+(r+2)^{-1}(r+1) (W_{r+1}-W_r)
\end{aligned} \mx{.}
\eeq
\end{lem}
The following notation will be used later
\beq{eq:notationz}
Z_r=:(T_2^r,Q_2^r,\ldots,T_N^r,Q_N^r) \mx{.}
\eeq
The assumptions $1.$ and $2.$ are nowadays ``standard'' and they have appeared in the literature, in various forms, since the earliest works. In particular (\ref{eq:hypdue}) is easily recognisable in  \cite{MR1375654} while (\ref{eq:hypuno}) can be classified as equivalent to the one appearing in \cite{MR1912262} and \cite{MR1996776} being related to the logarithmic time of transition. Assumption $3.$ is made for simplicity of discussion only.\\
In order to prove Lemma \ref{lem:convexity}, we need the following
\begin{prop}
Under the assumptions $1.$, $2.$ and $3.$ of Lemma \ref{lem:convexity}, recall (\ref{eq:tmenopiu}).  
Then, the following properties hold 
\begin{align}
\ml{T}^-(\mu) < 2 k_0 \K{k_0} & < \ml{T}^+(\mu) \label{eq:stimatmin}
\\
\Gamma_{i-1,i}^-(\omega_{i-1,i})+\Gamma_{i,i+1}^+(\omega_{i,i+1}) & =\ml{M}_{i-1,i+1}+o(1) \label{eq:boundsm} \mx{,}\\
\Theta_{i-1,i}^-(\omega_{i-1,i})+\Theta_{i,i+1}^+(\omega_{i,i+1}) & = o(1) \mx{,}
\label{eq:stimegammatheta} \\
(\pl_{T_2},\pl_{Q_2})  \Theta_{1,2}^+(\omega_{1,2})&=(- \Delta_{1,2}^{-1},1) \omega_{1,2} \drx{\Theta}_{1,2}^+ (\omega_{1,2}) +o(1) \mx{,} \label{eq:stimanuno}\\
(\pl_{T_N},\pl_{Q_N})  \Theta_{N,N+1}^-(\omega_{N,N+1})&=(\Delta_{N,N+1}^{-1}, -1) \omega_{N,N+1} \drx{\Theta}_{N,N+1}^-(\omega_{N,N+1})+ o(1) \mx{.}\label{eq:stimandue}
\end{align}
We have denoted $\ml{M}_{i,j}:=\ml{M}(\omega_{i,j})$ (understood $\omega_{i-1,i+1}:=(Q_{i+1}-Q_{i-1})/(T_{i+1}-T_{i-1})$), $\Delta_{i,j}:=T_{i}-T_{j}$ and, given $f=f(\omega,\ldots)$, $\drx{f}:=\pl_{\omega} f$.
\end{prop}
While (\ref{eq:stimatmin}) gives a lower and upper bounds for the time of transition as a function of $\mu$, (\ref{eq:boundsm}) and (\ref{eq:stimegammatheta}) reflect the assumption for which two consecutive frequencies are $\ml{O}(\mu)$ close, and that the instants $T_i$ are sufficiently large, as quantified by (\ref{eq:hypdue}) and (\ref{eq:hypuno}), respectively. In particular, as it is evident from the symmetry properties of the integrands in (\ref{eq:thetam}) and (\ref{eq:thetap}), one has  
\[
\lim_{(T_{i,i+1}^-,\omega_{i,i+1}) \rw (T_{i-1,i}^-,\omega_{i-1,i})} \left[\Theta_{i-1,i}^-(\omega_{i-1,i})+\Theta_{i,i+1}^+(\omega_{i,i+1})\right]=0 \mx{.} 
\]
\proof 
The estimate (\ref{eq:stimatmin}) follows from a result contained in \cite{carlson} and \cite{KARP2007186}. In particular, from  \cite[formulae (8) and (9)]{KARP2007186}, we get that 
\[
\K{k}=\frac{2}{\pi} \Kp{k} \log\frac{4}{k'}-\delta_1, \qquad \frac{(k')^2}{8} < \delta_1 < \frac{(k')^2 \log(4)}{k^2} \mx{.}
\]
First of all, from (\ref{eq:hypuno}) we get $0<\delta_1 < C \mu \log(4)$, hence
\[
2 k_0 \K{k_0} > \frac{2}{\sqrt{1+C \mu}} \left[ \frac{2}{\pi} \K{0}\log \left(4 \sqrt{\frac{1+C \mu}{C \mu}}\right)- C \mu \log(4) \right] \geq \frac{3}{2} \log \left(8 \sqrt{\frac{5}{\mu}}\right) \mx{.}
\]
The upper bound is found in a similar way, by using, in particular, $\K{\sqrt{C\mu/(1+C\mu)}} \leq 7/4$, as 
$\mu \leq \bar{\mu}$, by assumption. \\
As for (\ref{eq:boundsm}), let us observe that, by (\ref{eq:hypuno}), $(k_{i,i+1}')^2=\ml{O}(\mu)$. Hence, by (\ref{eq:resint}) and (\ref{eq:resto}), 
\begin{align}
\Gamma_{i-1,i}^-(\omega_{i-1,i})+\Gamma_{i,i+1}^+(\omega_{i,i+1}) & = 
\pi \left(\frac{\omega_{i-1,i}}{\sinh(\omega_{i-1,i} \pi/2)} +  
\frac{\omega_{i,i+1}}{\sinh(\omega_{i,i+1} \pi/2)}
\right) + o(1) \nonumber\\
& = \ml{M}(\omega_{i-1,i+1})+o(1)\nonumber
\end{align}
where (\ref{eq:hypdue}) has been used in the last passage.\\
The last bound, (\ref{eq:stimegammatheta}), is a straightforward check by using standard integral bounds. The only aspect worth pointing out is that, by recalling the expansion, see e.g. \cite[Pag. 574] {abramowitz1965handbook}
\beq{eq:expcn}
\cn(u,k)=\sech u + (k'/2)^2 \left[(u-\sinh u \cosh u)\tanh u \sech u \right]+ \ml{O}\left((k')^4\right) \mx{,}
\eeq
the smallness of the difference $\cn(u,k)-\sech u$, as $k \rw 1$ (in some suitable norm), follows from the analyticity of $\cn(u,k)$, see \cite{10.2307/3560143}. More precisely, one can bound, for all $u:=k^{-1}t \in [-\K{k},\K{k}]$,  
\[
\left|(k'/2)^2 \left[(u-\sinh u \cosh u)\tanh u \sech u \right] \right| \leq 4^{-1} (k')^2 e^{-u} \leq 5 \sqrt{\mu/2} \mx{,}
\]
where (\ref{eq:hypuno}) and (\ref{eq:stimatmin}) have been used in the last inequality.\\
We shall only show how to prove (\ref{eq:stimanuno}) for $\pl_{T_2}$, the rest of the proofs being analogous. By differentiating (\ref{eq:k}) one gets
\beq{eq:diffk}
\pl_{T_i} \left({k_{i-1,i},k_{i,i+1}} \right)=2^{-1} \left(\E{k_{i-1,i}}^{-1}(1-k_{i-1,i}^2), -\E{k_{i,i+1}}^{-1}(1-k_{i,i+1}^2) \right) \mx{,}
\eeq
see \cite[Pag. 521]{whittaker1996course}. The latter quantities, by (\ref{eq:hypuno}), are $\ml{O}(\mu)$, hence
\begin{align}
\pl_{T_2} \Theta_{1,2}^+(\omega_{1,2}) &=\cn^2(\K{k_{1,2}},k_{1,2})\sin(\omega_{1,2}k_{1,2}\K{k_{1,2}}) \nonumber\\
&+ 2 \pl_{T_2}k_{1,2} \int_0^{T_{1,2}^-}\pl_{k_{1,2}} \left[\cn^2(k_{1,2}^{-1}t,k_{1,2})  \sin(\omega_{1,2} t) \right]dt \nonumber \\
&+ 2 \pl_{T_2} \omega_{1,2}\int_0^{T_{1,2}^-} \cn^2(k_{1,2}^{-1}t,k_{1,2}) t \cos(\omega_{1,2} t) dt \mx{.}
\end{align}
The first term is zero as $\cn(\K{k},k)=0$ and the second one is $o(1)$ by (\ref{eq:diffk}) and the analyticity of $\cn(k^{-1} t,k)$ for $k \in (k_0,1)$. The last integral is, by definition, $\drx{\Theta}_{1,2}^+ (\omega_{1,2})$. The proof is now complete.
\endproof 
\noindent 
We now recall a standard result concerning Nesterov's method by a straightforward adaptation of a statement contained in \cite{JMLR:v17:15-084}. It should be mentioned that \cite{JMLR:v17:15-084} shows a remarkable interpretation of such a method in terms of the discrete version of a suitable differential equation. 
\begin{prop}[Nesterov]\label{prop:nesterov}
Let $\ml{U}\subset \RR^n$ be an open set and let $\ml{F}_n:\ml{U} \rw \RR$ be a convex Lipschitz function with Lipschitz constant $L$. Suppose that there exists $\hat{Z} \in \ml{U}$ such that $\nabla \ml{F}_n (\hat{Z})=0$. Let us denote $\hat{\ml{F}}_n:=\ml{F}_n (\hat{Z})$.\\
Then, for all $Z_0 \in \ml{U}$ and all $\alpha \leq 1/L$, the scheme (\ref{eq:graddesc}) (in which $\ml{F}$ is replaced by $\ml{F}_n$) converges to $\hat{Z}$, which is the unique minimum of $\ml{F}_n$ over $\ml{U}$. Furthermore, the rate of convergence is bounded as follows
\beq{eq:speedofnesterov}
\ml{F}_n(Z_r)-\hat{\ml{F}}_n \leq (2/\alpha)(r+1)^{-2}|Z_r-Z_0|^2 \mx{.}
\eeq 
\end{prop}
\proof [Proof of Lemma \ref{lem:convexity}]
The strategy consists in showing that on $\mathfrak{C}_{\ml{I},\pi/4}$:
\begin{enumerate}
\item Step $1$: The gradient of $\ml{F}^{(\leq 1)}$ satisfies the assumption of the Poincar\'{e}-Miranda (PM) Theorem \cite{mira}, see also \cite{kulpa}; and 
\item Step $2$: $\ml{F}^{(\leq 1)}$ is convex. 
\end{enumerate}
Defined ${C}_l:=[-l,l]^n$ with $l>0$, and given $f_i:{C}_l \rw \RR$, $i=1,\ldots,n$ continuous functions such that $f_j(x)|_{x_j=-l} \leq 0$ and $f_j(x)|_{x_j=l} \geq 0$ for all $j=1,\ldots,n$, let us briefly recall that the PM Theorem establishes the existence of at least one point $x^* \in C$ such that $f_i(x^*)=0$ for all $i$. Furthermore, it is immediate that $x^* \in C_l\setminus \partial C_l$ if both of the above mentioned conditions are satisfied with strict inequalities.\\
We will apply this result to the gradient of $\ml{F}^{(\leq 1)}$ showing the existence of \emph{at least} one stationary point. Afterwards, Prop. \ref{prop:nesterov} on $\mathcal{U}:=\mathfrak{C}_{\ml{I},\pi/4} \setminus \pl \mathfrak{C}_{\ml{I},\pi/4}$ will be used to conclude the proof on $\ml{F}$, for sufficiently small $\mu$.
\subsubsection*{Step $1$}
Firstly, we observe that (\ref{eq:stimegammatheta}) implies $\Gamma_{i-1,1}^-(1)+\Gamma_{i,i+1}^+(1)=\ml{M}(1)+o(1)$ then we set $\ml{A}:=\ml{M}(1)$. Hence, by (\ref{eq:boundsm}) and (\ref{eq:stimegammatheta}) we get
\begin{align}
\ml{F}^{(1)}&= \Theta_{1,2}^+(\omega_{1,2})\sin Q_1+ \Theta_{N,N+1}^-(\omega_{N,N+1})\sin Q_{N+1} \nonumber \\
&-2^{-1} \left[\ml{M}_{1,2} \cos Q_1+\ml{M}_{N,N+1} \cos Q_{N+1} \right] \nonumber \\
&-\sum_{i=2}^N \left[ \ml{A}\cos T_i + \ml{M}_{i-1,i+1}\cos Q_i\right]+\ml{C}+o(1) \label{eq:lastfuno} 
\end{align}
where
\[
\ml{C}:=\Theta_{1,2}^+(1)\sin T_1-\Gamma_{1,2}^+(1)\cos T_1+
\Theta_{N,N+1}^-(1)\sin T_{N+1}-\Gamma_{N,N+1}^-(1)\cos T_{N+1} \mx{.}
\]
Roughly speaking, the convexity of $\ml{F}^{(\leq 1)}$ relies on the terms $\ml{A}\cos T_i + \ml{M}_{i-1,i+1}\cos Q_i$ and the fact that the objects $\ml{A}$ and $\ml{M}_{i-1,i+1}$ are bounded by below, reflecting the contribution of the Melnikov term.\\
From the proof of (\ref{eq:stimanuno}), it is easy to show that $\pl_{T_2}  \Theta_{1,2}^+(1)=o(1)$ and $\pl_{T_N}  \Theta_{N,N+1}^-(1)=o(1)$. Hence, by using these, together with (\ref{eq:stimanuno}) and (\ref{eq:stimandue}), one gets, up to $\ml{O}(\mu)$
\begin{align*}
\pl_{T_2} \ml{F}^{(1)} & =\ml{A} \sin T_2- \Delta_{2,4}^{-1} \omega_{2,4} \drx{\ml{M}}_{2,4} \cos Q_3 + (2 \Delta_{1,2})^{-1} \omega_{1,2} \drx{\ml{M}}_{1,2} \cos Q_1 \\
&-\Delta_{1,2}^{-1} \omega_{1,2} \drx{\Theta}_{1,2}^+(\omega_{1,2}) \sin Q_1 \\
\pl_{Q_2} \ml{F}^{(1)} & =\ml{M}_{1,3}\sin Q_2 + \Delta_{2,4}^{-1}  \drx{\ml{M}}_{2,4} \cos Q_3 - (2 \Delta_{1,2})^{-1} \omega_{1,2} \drx{\ml{M}}_{1,2} \cos Q_1 \\
&+\Delta_{1,2}^{-1}\drx{\Theta}_{1,2}^+(\omega_{1,2}) \sin Q_1
\\
& \dots \\
\pl_{T_i}\ml{F}^{(1)}& =\ml{A} \sin T_i - \Delta_{i,i+2}^{-1} \omega_{i,i+2} \drx{\ml{M}}_{i,i+2} \cos Q_{i+1}+ \Delta_{i-2,i}^{-1} \omega_{i-2,i} \drx{\ml{M}}_{i-2,i} \cos Q_{i-1}\\
\pl_{Q_i}\ml{F}^{(1)}& =
\ml{M}_{i-1,i+1} \sin Q_i - \Delta_{i,i+2}^{-1} \omega_{i,i+2} \drx{\ml{M}}_{i,i+2} \cos Q_{i+1}+ \Delta_{i-2,i}^{-1} \omega_{i-2,i} \drx{\ml{M}}_{i-2,i} \cos Q_{i-1}\\
& \ldots \\
\pl_{T_N}\ml{F}^{(1)}& =\ml{A}\sin T_N - \Delta_{N-2,N}^{-1}\omega_{N-2,N} \drx{\ml{M}}_{N-2,N}\cos Q_{N-1}\\
& + (2 \Delta_{N,N+1})^{-1} \omega_{N,N+1}\drx{\ml{M}}_{N,N+1} \cos Q_{N+1}\\&+\Delta_{N,N+1}^{-1} \omega_{N,N+1} \drx{\Theta}_{N,N+1}^-(\omega_{N,N+1}) \sin Q_{N+1} \\
\pl_{Q_N}\ml{F}^{(1)}& = \ml{M}_{N-1,N+1} \sin Q_N + (2 \Delta_{N,N+1})^{-1}\drx{\ml{M}}_{N,N+1}\cos Q_{N+1} \\
&-\Delta_{N,N+1}^{-1} \drx{\ml{M}}_{N-2,N}\cos Q_{N-1}-\Delta_{N,N+1}^{-1}  \drx{\Theta}_{N,N+1}^-(\omega_{N,N+1}) \sin Q_{N+1} 
\end{align*}
where $i=2,\ldots,N-1$. \\
On the other hand, it is immediate from (\ref{eq:fzero}) and the elementary derivation rules for elliptic integrals, that 
\begin{align*}
\pl_{T_i} \ml{F}^{(0)}& = 2\left(k_{i,i+1}^{-1}-k_{i-1,i}^{-1} \right) + 2^{-1} \left( \omega_{i,i-1}^2-\omega_{i,i+1}^2 \right)\\
\pl_{Q_i} \ml{F}^{(0)} & = \omega_{i+1,i}-\omega_{i,i+1}
\end{align*}
for all $i=2,\ldots,N$, which compares straightforwardly with the expression for the ``full'' $\mathcal{F}$ given by (\ref{eq:gradbertiuno}) and (\ref{eq:gradbertidue}). \\
The latter clearly shows the correspondence between the stationary points of $\mathcal{F}$ and ``true'' solutions for the system (\ref{eq:lag}) i.e. whose speeds are continuous also at the instants $t=T_i$. This property has been explicitly used since \cite{MR1912262}.   
The verification of the hypotheses of the PM Theorem is similar for all the components of the gradient. We sketch how to proceed for $\pl_{T_2} \ml{F}$ for completeness. \\
First of all we notice that (\ref{eq:melnikov}) and (\ref{eq:hyptre}) implies
\beq{eq:stimem}
\left| \drx{\ml{M}} \right| \leq 4, \qquad \left| \ddrx{\ml{M}} \right| \leq 7/2 \mx{.}
\eeq
Secondly, by recalling (\ref{eq:expcn}), one has  
\[
\left|\drx{\Theta}_{1,2}^+(\omega_{1,2})\right| \leq k \int_{0}^{+\infty}  u \sech^2 u du \leq 1, \qquad  \left|\ddrx{\Theta}_{1,2}^+(\omega_{1,2})\right| \leq 4 \mx{.}
\]  
Similarly for the derivatives of $\Theta_{N,N+1}^-$. Hence 
\[
\left|\Delta_{2,4}^{-1} \omega_{2,4} \drx{\ml{M}}_{2,4} \cos Q_3 + (2 \Delta_{1,2})^{-1} \omega_{1,2} \drx{\ml{M}}_{1,2} \cos Q_1
-\Delta_{1,2}^{-1} \omega_{1,2} \drx{\Theta}_{1,2}^+(\omega_{1,2}) \sin Q_1 \right| \leq 5/\Delta_- \mx{,}
\]
where $\Delta_-:=\min_{i=2,\ldots,N} (T_{i}-T_{i-1})$. By (\ref{eq:stimatmin}), the latter satisfies $\Delta_- > 47/4$, being $\mu \leq \bar{\mu}$.\\
On the other hand, by (\ref{eq:hypuno}), (\ref{eq:hypdue}) and the chosen value of $C$, we have $\left|\ml{F}^{(0)}\right|\leq \mu/8 $. Therefore, 
\beq{eq:previousbound}
\pl_{T_2} \ml{F}^{(\leq 1)}|_{T_2=-\pi/4}=\mu\left[-\ml{A}\sin \left(\pi/4 \right) + 1/8 +20/47\right]<-\mu<0 \mx{,}
\eeq
for all the involved $T_i,Q_i \in [-\pi/4,\pi/4]$, (recall $\ml{A}:=\ml{M}(1)=2 \pi/ \sinh (\pi/2) \sim 2.73$). Analogously, $\pl_{T_2} \ml{F}^{(\leq 1)}|_{T_2=\pi/4}>0$. The procedure for the other variables is similar.\\
Hence, by the PM Theorem, the function $\ml{F}$ possesses at least one stationary point on $\mathfrak{C}_{\ml{I},\pi/4} \setminus \pl \mathfrak{C}_{\ml{I},\pi/4}$, for sufficiently small $\mu$.\\
As it will be used later, we notice that the bound 
\beq{eq:boundnablaf}
\sup_{Z \in \mathfrak{C}_{\ml{I},\pi/4}} \enorm{\nabla \ml{F}^{(\leq 1)}(Z)} \leq 4 \mu \mx{,}
\eeq
can be obtained exactly in the same way as (\ref{eq:previousbound}). 
\subsubsection*{Step 2}
First of all let us notice that the set 
$\mathfrak{C}_{\ml{I},\pi/4}$ is convex (straightforward check), then we split $\ml{F}^{(\leq 1)}$ as follows
\beq{eq:split}
\ml{F}^{(\leq 1)}=\Omega+W, \qquad 
\eeq
where
\begin{align*}
\Omega&:=2 \sum_{i=1}^N \omega_{i,i+1}^2 T_{i,i+1}^- \\
W&:=(\ml{F}^{(0)}-\Omega)+\ml{F}^{(1)}
\end{align*}
with $\ml{F}^{(0)}$ and $\ml{F}^{(1)}$ of the form (\ref{eq:fzero}) and  (\ref{eq:lastfuno}), respectively.\\
We shall proceed by showing that, under the current assumptions, $\Omega$ is convex and $W$ is strictly convex. In such a way, $\ml{F}^{(\leq 1)}$ is strictly convex as well. The splitting (\ref{eq:split}) is motivated by the fact that, as we shall see in a moment, the second derivatives of $\Omega$ are $\ml{O}(1)$ while those of $W$ are $\ml{O}(\mu)$.\\
Let us start from the strict convexity of $W$. In order to show this property, it is sufficient to check that its Hessian matrix, $\pl^2 W$, has a strictly dominant diagonal. \\
By symmetry, and obvious regularity properties, the matrix $\{\pl^2 W\}_{l,m}$ is defined once the elements with $m \geq l$ have been written. In particular, the Hessian matrix of $\ml{F}^{(0)}-\Omega$ is entirely defined by the following two elements, as $m=2,\ldots,N-1$,
\begin{align*}
\pl_{T_m^2}^2 (\ml{F}^{(0)}-\Omega)&=\frac{1-k_{m-1,m}^2}{k_{m-1,m}^3 \E{k_{m-1,m}}}-\frac{1-k_{m,m+1}^2}{k_{m,m+1}^3 \E{k_{m,m+1}}}\\ 
\pl_{T_m T_{m+1}}^2 (\ml{F}^{(0)}-\Omega)&=-\frac{1-k_{m,m+1}^2}{k_{m,m+1}^3 \E{k_{m,m+1}}} 
\end{align*}
which are clearly $\ml{O}(\mu)$ due to (\ref{eq:hypuno}).\\
On the other hand, the structure of $\pl^2 \ml{F}^{(1)}$ is more involved. More precisely, we have: 
\begin{align*}
\pl_{T_2^2}^2 \ml{F}^{(1)} & = \ml{A}\cos T_2 +\Delta_{1,2}^{-2}\omega_{1,2}(\omega_{1,2} \ddrx{\Theta}_{1,2}^-+2\drx{\Theta}_{1,2}^-)\sin Q_1
\\
& - (2 \Delta_{1,2}^{2})^{-1} \omega_{1,2}(\omega_{1,2} \ddrx{\ml{M}}_{1,2}+2\drx{\ml{M}}_{1,2})\cos Q_1
+ \Delta_{2,4}^{-2}\omega_{2,4}(\omega_{2,4} \ddrx{\ml{M}}_{2,4}+2\drx{\Theta}_{2,4})\cos Q_3 \\
\pl_{T_2 Q_2}^2 \ml{F}^{(1)} & = -\Delta_{1,2}^{-2}(\omega_{1,2} \ddrx{\Theta}_{1,2}^-+\drx{\Theta}_{1,2}^-)\sin Q_1 +
(2 \Delta_{1,2}^{2})^{-1} (\omega_{1,2} \ddrx{\ml{M}}_{1,2}+\drx{\ml{M}}_{1,2})\cos Q_1\\
& + \Delta_{2,4}^{-2}(\omega_{2,4} \ddrx{\ml{M}}_{2,4}+2\drx{\Theta}_{2,4}^-)\cos Q_3\\
\pl_{Q_2^2}^2 \ml{F}^{(1)} & = \ml{M}_{1,3}\cos Q_2 +\Delta_{1,2}^{-2}\drx{\Theta}_{1,2}^- \sin Q_1
- (2 \Delta_{1,2}^{2})^{-1} \ddrx{\ml{M}}_{1,2}\cos Q_1
+ \Delta_{2,4}^{-2} \ddrx{\ml{M}}_{2,4} \cos Q_3 
\end{align*}
then for all $m=2,\ldots,N-2$ 
\begin{align*}
\pl_{T_m Q_{m+1}}^2 \ml{F}^{(1)} & = \Delta_{m,m+2}^{-1} \omega_{m,m+2} \drx{\ml{M}}_{m,m+2} \sin Q_{m+1} \\
\pl_{T_m T_{m+2}}^2 \ml{F}^{(1)} & = \Delta_{m,m+2}^{-2} \omega_{m,m+2}
(\omega_{m,m+2} \ddrx{\ml{M}}_{m,m+2}+2\drx{\ml{M}}_{m,m+2})\cos Q_{m+1}\\
\pl_{T_m Q_{m+2}}^2 \ml{F}^{(1)} & =-\Delta_{m,m+2}^{-2} 
(\omega_{m,m+2} \ddrx{\ml{M}}_{m,m+2} +\drx{\ml{M}}_{m,m+2})\cos Q_{m+1}\\
\pl_{Q_m T_{m+1}}^2 \ml{F}^{(1)} & = -\Delta_{m_1,m+1}^{-1} \omega_{m-1,m+1} \drx{\ml{M}}_{m-1,m+1}\sin Q_m\\
\pl_{Q_m Q_{m+1}}^2 \ml{F}^{(1)} & = \Delta_{m-1,m+1}^{-1} \drx{M}_{m-1,m+1} \sin Q_m - \Delta_{m,m+2}^{-1} \drx{M}_{m,m+2} \sin Q_{m+2} \\
\pl_{Q_m T_{m+2}}^2 \ml{F}^{(1)} & = -\Delta_{m,m+2}^{-2} 
(\omega_{m,m+2} \ddrx{\ml{M}}_{m,m+2} +\drx{\ml{M}}_{m,m+2})\cos Q_{m+1}\\
\pl_{T_m Q_{m+2}}^2 \ml{F}^{(1)} & =\Delta_{m,m+2}^{-2} 
 \ddrx{\ml{M}}_{m,m+2} \cos Q_{m+1}
\end{align*}
and all $m=3,\ldots,N-2$. 
\begin{align*}
\pl_{T_m^2}^2 \ml{F}^{(1)} & = \ml{A}\cos T_m -
 \Delta_{m-2,m}^{-2} \omega_{m-2,m}(\omega_{m-2,m} \ddrx{\ml{M}}_{m-2,m}+2\drx{\ml{M}}_{m-2,m})\cos Q_{m-1}
\\
&- \Delta_{m,m+2}^{-2}\omega_{m,m+2}(\omega_{m,m+2} \ddrx{\ml{M}}_{m,m+2}+2\drx{\ml{M}}_{m,m+2})\cos Q_{m+1} \\
\pl_{T_m Q_m}^2 \ml{F}^{(1)} & = \Delta_{m-2,m}^{-2}(\omega_{m-2,m} \ddrx{\ml{M}}_{m-2,m}+\drx{\ml{M}}_{m-2,m})\cos Q_{m-1} \\
& + \Delta_{m,m+2}^{-2} (\omega_{m,m+2} \ddrx{\ml{M}}_{m,m+2}+\drx{\ml{M}}_{m,m+2})\cos Q_{m+1}\\
\pl_{Q_m^2}^2 \ml{F}^{(1)} & = \ml{M}_{m-1,m+1}\cos Q_m 
- \Delta_{m-2,m}^{-2} \ddrx{\ml{M}}_{m-2,m}\cos Q_{m-1}
- \Delta_{m,m+2}^{-2} \ddrx{\ml{M}}_{m,m+2}\cos Q_{m+1} \\
\end{align*}
The entries for $m=N-1$ and $m=N$ are straightforward, while the remaining ones are either zero or obtained by symmetry.\\
The diagonal dominance of the decadiagonal Hessian matrix of $W$ is proven once we show that the off-diagonal entries are small enough with respect to the the positive contribution of the terms $\ml{A}\cos T_m$ or $\ml{M}_{m-1,m+1}\cos Q_m$. For this purpose, a suitable ``smallness'' of the terms $\Delta_{m-1,m+1}^{-1}$, $\sin T_m$, $\sin Q_m$ and $C^{-1}$ (the latter concerning $\ml{F}^{(0)}-\Omega$), plays a key role. An example of such an estimate is given below for completeness, the others are similar. \\
Let us deal with the sixth row of $\pl^2 \ml{F}^{(1)}$ (which is the first full row from the top). First of all, by (\ref{eq:stimatmin}), we get $\min\{\Delta_{4,6},\Delta_{5,7},\Delta_{6,8}\} \geq 4 k_0 \K{k_0} \geq 47/2$ for all $\mu \leq \bar{\mu}$. Hence, taking into account of the symmetry of $\pl^2 \ml{F}^{(1)}$, then using the bounds (\ref{eq:stimem}), the sum of the the off-diagonal terms of the sixth row of $\pl^2\ml{F}^{(1)}$ can be bounded as follows 
\beq{eq:sixthone}
\sum_{l=1}^{10} \left| \pl_{Q_4 Z_l}^2 \ml{F}^{(1)}\right| \leq \frac{150}{2209} + \frac{42}{47 \sqrt{2}} \mx{,}
\eeq
where $Z:=(T_2,Q_2,\ldots,T_6,Q_6)$. On the other hand, by recalling the value of $C$, a bound on the off-diagonal terms of $\ml{F}^{(0)}-\Omega$, on the same row, is 
\beq{eq:sixthtwo}
\left|\frac{1-k_{6,7}^2}{k_{6,7}^3 \E{k_{6,7}}} \right|+ \left|\frac{1-k_{7,8}^2}{k_{7,8}^3 \E{k_{7,8}}} \right| \leq \frac{1}{10} \mx{.}
\eeq
By considering the negative terms on the diagonal of both matrices and using (\ref{eq:sixthone}) and (\ref{eq:sixthtwo}), a sufficient condition for the strict diagonal dominance is 
\beq{eq;conconvw}
\ml{A}\cos Q_6 - 28/2209- 1/20>  1/10+150/2209+42/(47 \sqrt{2}) \mx{,}
\eeq
which is satisfied as $Q_6 \Mod 2 \pi \in [0,\pi/4]\cup[3\pi/4,2 \pi)$. This proves the strict convexity of $W$. \\
We can now proceed by showing the convexity of $\Omega$. For this purpose, let us suppose to perform the linear transformation of variables given by the following formulae
\[
\eta_l:=T_{l}-T_{l-1}, \qquad \xi_l:=Q_{l}-Q_{l-1} \mx{,}
\]
for all $l=3,\ldots,N$ and $(\eta_2,\xi_2):=(T_2,Q_2)$. This is meaningful as convexity is invariant (in particular) by linear transformations of variables, see e.g. \cite{convex}. \\
The determinant of the matrix associated to the transformation is equal to one, in particular the transformation is invertible and we obtain, in the new set of variables
\[
\Omega=\frac{1}{2} \left[\frac{\eta_2^2}{\xi_2} + \frac{\eta_3^2}{\xi_3} + \frac{\left( Q_{N+1}-\sum_{j=2}^N \xi_j\right)^2}{\left( T_{N+1}-\sum_{j=2}^N \eta_j \right)}\right] \mx{,}
\]  
where $Q_{N+1}$, $T_{N+1}$ are constants and $Q_1=T_1=0$ without loss of generality. The proof is complete as the function of two variables $f(x,y):=y^2/x$ is convex for all $x \neq 0$, see \cite[P.73]{convex}, so that $\Omega$ is convex, being the sum of convex functions.\\
The last task consists in the computation of a bound for $L$. It is important to notice that, whilst $\pl^2 W$ contains all $O(\mu)$ terms, the quantity $\norm{\pl^2 \Omega}$ is not smaller than $O(\log \mu^{-1})$ due to terms bounded by $1/\Delta_-$ arising in $\pl^2 \Omega$ (straightforward computation). This implies that, for ``realistic'' values of $\mu$ (i.e. of the same magnitude of $\bar{\mu}$), the constant $L$ is not ``that small'' and a refined bound would be unnecessary for our purposes. Checking that one can choose $L = 1$ is similar to the bounds above for $\pl^2 W$ and it is just a matter of patience. This completes the proof.     
\endproof

\section{The Boundary Value Problem}\label{sec:bvp}
In this section we deal with the problem (\ref{eq:bvps}) in its discrete form, i.e. with the minimisation of the function $\ml{F}_{i,i+1}|_{\ml{M}_{i,i+1}}$. It is immediate to realise that (\ref{eq:bvpuno}) (e.g. with $l=1$) is a prototype for (\ref{eq:bvps}) as $i=1,\ldots,N$, hence the whole analysis can be reduced to to the study of (\ref{eq:bvpuno}). Setting $\ml{F}_{a,b}:=\int_{T_a}^{T_b} \ml{L}_\mu$ and 
\[
\ml{M}_{a,b}:=\{t_k=T_a+\tilde{h} k,\, k=1,\ldots,\tilde{n} \}, \qquad (\tilde{n}+1) \tilde{h}:=\Delta_{a,b} \mx{,} 
\]
where $\Delta_{a,b}:=T_b-T_a$, we are interested in finding the minima $(\underline{q},\underline{Q})$ of $\ml{F}_{a,b}|_{\ml{M}_{a,b}}$. \\
Let us set 
\beq{eq:settdiscr}
\ue{q}=:\ue{q}^0+\ue{v},\qquad \ue{Q}=:\ue{Q}^0+\ue{w}\mx{,}
\eeq 
where $(\ue{q}^0,\ue{Q}^0)=(q_0(t),Q_0(t))|_{\ml{M}_{a,b}}$. Note that the form (\ref{eq:settdiscr}) turns out to be particularly convenient from the computational point of view, as either $q_0(t)$ or $Q_0(t)$ are explicitly known. In particular, $q_0(t)$ can be approximated via dedicated in-built packages for Jacobian elliptic functions, see e.g. \cite{oct} or known approximation formulae based on Fourier expansions, see \cite[Pag. 510]{whittaker1996course}. \\
Recalling (\ref{eq:discretef}) and (\ref{eq:discretel}), the stationary points of $(\underline{q},\underline{Q})$ of $\ml{F}_{a,b}|_{\ml{M}_{a,b}}$ must be solutions of the following system of $2 \tilde{n}$ non-linear equations
\beq{eq:discreteel}
\underline{\Psi}(\underline{v},\underline{w},\ue{t};\mu):=
\begin{pmatrix}
	\underline{\underline{T}} & 0\\
	0 & \underline{\underline{T}}\\
\end{pmatrix}
\begin{pmatrix}
	\underline{v} \\
	\underline{w} 
\end{pmatrix}
+\tilde{h}^2
\begin{pmatrix}
	\underline{G} + \mu \underline{V}\\
	\mu \underline{W}
\end{pmatrix}=\underline{0}
\mx{,}
\eeq
where $\underline{\underline{T}}:=\trid(-1,2,-1) \in GL(\tilde{n},\ZZ)$, and for all $k=1,\ldots,\tilde{n}$,
\begin{align*}
G_k&:=\sin (q_{k}^0+v_k)-\sin q_k^0 \\
V_k&:=-\sin (q_k^0+v_k) \left[ \cos (Q_k^0+w_k)+\cos t_k \right]  \\
W_k&:=\sin (Q_k^0+w_k) \left[1- \cos (q_k^0+v_k) \right] \mx{.}
\end{align*}
It is a standard result, see \cite{mawe}, that the system (\ref{eq:discreteel}) is nothing but the \emph{discrete Euler-Lagrange equations} obtained via the standard discretisation of the second derivative in  
\beq{eq:lageqvw}
\begin{cases}
	\ddot{v}&=\sin (q_0+v)-\sin q_0-\mu \sin (q_0+v) [\cos (Q_0+w)+\cos t]\\
	\ddot{w}&=\mu \sin (Q_0+w) [1-\cos (q_0+v)]
\end{cases}
\mx{,}
\eeq
subject to $v(T_a)=w(T_a)=v(T_b)=w(T_b)=0$, i.e. from equations associated with (\ref{eq:lag}) (under the known boundary conditions) by setting $q(t)=q_0(t)+v(t;\mu)$ and $Q(t)=Q_0(t)+w(t;\mu)$.\\
The Jacobian matrix of $\ue{\Psi}$ can be written as 
\beq{eq:jacobianpsi}
\frac{\pl}{\pl(\ue{v},\ue{w})} \ue{\Psi} 
=
\begin{pmatrix}
	\underline{\underline{T}} & \ue{\ue{0}}\\
	\ue{\ue{0}} & \underline{\underline{T}}\\
\end{pmatrix}
+\tilde{h}^2
\left\{
\begin{pmatrix}
	\underline{\underline{D}}(\ue{v},t) & \ue{\ue{0}}\\
	\ue{\ue{0}} & \underline{\underline{0}}\\
\end{pmatrix}
+
\mu
\begin{pmatrix}
	\underline{\underline{A}} & \ue{\ue{B}}\\
	\ue{\ue{B}} & \underline{\underline{C}}\\
\end{pmatrix} (\ue{v},\ue{w},\ue{t})
\right\}
=:\ue{\ue{\ml{T}}}+\tilde{h}^2 \left\{ \ue{\ue{\ml{D}}}+\mu \ue{\ue{\ml{V}}}\right\}
\eeq
where $\ue{\ue{A}}, \ue{\ue{B}}, \ue{\ue{C}}$ and $\ue{\ue{D}}$ are diagonal matrices, whose diagonal elements read as
\begin{align*}
a_{k,k}&:=\cos (q_{k}^0+v_k) \left[ \cos (Q_k^0+w_k)+\cos t_k \right]\\
b_{k,k}&:=\sin (q_k^0+v_k) \sin (Q_k^0+w_k)  \\
c_{k,k}&:=\cos (Q_k^0+w_k) \left[1- \cos (q_k^0+v_k) \right]\\
d_{k,k}&:=\cos (q_{k}^0+v_k)
\end{align*}
respectively. \\
As it is easy to check, we have $\ue{\Psi}(\ue{0},\ue{0},\ue{t};0)=\ue{0}$. Hence, it is natural to ask whether or not the Jacobian matrix of $\ue{\Psi}$ computed at this point, i.e.  
\beq{eq:jacobianpsizero}
\ue{\ue{\mathcal{J}_0}}:=\frac{\pl}{\pl(\ue{v},\ue{w})} \ue{\Psi} (\ue{0},\ue{0},\ue{t};0)=
\ue{\ue{\ml{T}}}+\tilde{h}^2 \ue{\ue{\ml{D}}}(\ue{0},\ue{t})
\mx{,}
\eeq
is non-singular, allowing to use, in this way, the Implicit Function Theorem (IFT). The answer is (more than) affirmative, as stated in the next
\begin{prop}\label{prop:defpos}%
Suppose that $\Delta_{a,b} \geq 3 \pi$ and that $\tilde{h} \leq 0.01$ (via a suitable choice of $\tilde{n}$), then 
\beq{eq:lowesteig}
\lambda^{-}(\ue{\ue{\mathcal{J}_0}}) \geq (\tilde{\alpha}/\Delta_{a,b})^2 \tilde{h}^2\mx{,}
\eeq
where $\tilde{\alpha}:=3 \pi /4$. In particular,  $\ue{\ue{\mathcal{J}_0}}$ is positive definite.
\end{prop}
\proof
First of all we recall the classical result, see e.g. \cite{meyer2000matrix}, which states that the eigenvalues of the matrix $\ue{\ue{T}} \in GL(\tilde{n},\ZZ)$ read as $
\lambda_j=4 \sin^2 [(2 \tilde{n} + 2)^{-1} j \pi ]$, for all $j=1,\ldots,\tilde{n}$. A possible bound for the smallest of them is 
\beq{eq:smallesteigt}
\lambda^-(\ue{\ue{T}}) \equiv \lambda_1 \geq 9 \pi^2[4(\tilde{n}+1)]^{-2} \mx{,}
\eeq
implying the (well known) positive definiteness of $\ue{\ue{T}}$. On the other hand, due to the structure of $\ue{\ue{\mathcal{J}_0}}$, supposing $\lambda^-(\ue{\ue{T}}+\tilde{h}^2 \ue{\ue{D}}(0,\ue{t}))>0$, we have $\lambda^-(\ue{\ue{\mathcal{J}_0}}) = \min\{\lambda^-(\ue{\ue{T}}),\lambda^-(\ue{\ue{T}}+\tilde{h}^2 \ue{\ue{D}}(0,\ue{t}))\}$. Hence, hence we need to establish whether or not $\lambda^-(\ue{\ue{T}}+\tilde{h}^2 \ue{\ue{D}}(0,\ue{t}))$ is strictly positive and, if this is the case, find an estimate from below. The difficulty in doing so is suggested by the Gershgorin Theorem and subsequent generalisations: as $q_k^0$ varies between $\pi$ and $3 \pi$ there are diagonal elements of $(\ue{\ue{T}}+\tilde{h}^2 \ue{\ue{D}}(0,\ue{t}))$, i.e. $2+\tilde{h}^2\cos q_{k}^0$, that are strictly smaller than $2$ and there is no obstruction for a set of eigenvalues to be located in the negative (real) semi-axis. The key point consists showing that the elements such that $\cos q_{k}^0<0$ are ``not too many'' with respect to those for which $\cos q_{k}^0>0$. Roughly speaking, this can be shown to be true as the pendulum ``spends most of the time close to the unstable equilibrium'' $2 \pi$. The floating point notation is used in this proof (where convenient) for the sake of simplicity, e.g. $0.52$ will be used instead of $52/100$. \\
First of all, let us recall the well known \emph{Sturm sequence property}, see e.g. \cite{BUCHHOLZER20121837} for a symmetric tridiagonal matrix $\ue{\ue{\mathfrak{T}}}:=\trid(\ue{\beta},\ue{\alpha},\ue{\beta}) \in GL(n,\RR)$ where $\ue{\alpha}=(\alpha_1,\ldots,\alpha_n)$ and $\ue{\beta}=(\beta_1,\ldots,\beta_{n-1})$. In particular, by setting 
\[
\begin{cases}
p_r(x)& = (\alpha_r-x)p_{r-1}(x)-\beta_{r-1}^2 p_{r-2}(x),\qquad r=1,\ldots,n\\
p_0(x) & \equiv 1
\end{cases}
\mx{,}
\]
it is possible to show that the number of sign changes of the sequence $\{p_m\}_{m=0,\ldots,n}$ equals the number of eigenvalues of $\ue{\ue{\mathfrak{T}}}$ smaller than $x$. \\
Clearly, $(\ue{\ue{T}}+\tilde{h}^2 \ue{\ue{D}}(0,\ue{t}))$ is in the described form as it equals $\ue{\ue{\mathfrak{T}}}$ with $\alpha_r:=d_{r,r}|_{v_r=0}$ and $\beta_r \equiv -1$ for $n=\tilde{n}$.
Let us suppose we want to show that $\ue{\ue{\mathfrak{T}}}$ has eigenvalues not smaller than $\tilde{h}^2/16$. If this is true, recalling the definition of $\tilde{h}$ and (\ref{eq:smallesteigt}), under the assumption on $\Delta_{a,b}$, we have that $\lambda^-(\ue{\ue{T}}) \leq \lambda^-(\ue{\ue{T}}+\tilde{h}^2 \ue{\ue{D}}(0,\ue{t}))$, which gives (\ref{eq:lowesteig}).\\ According to the mentioned result on the Sturm sequence property, it is sufficient to prove that the sequence $\{q_r\}_{r=1,\ldots,\tilde{n}}$ with $q_r:=p_r(\tilde{h}^2/16)$ is non-negative. For this purpose, let us define $k_{a,b}$ as satisfying $\Delta_{a,b}=2 k_{a,b} \K{k_{a,b}}$, then recall that it is possible to perform a sequence of (similarity) transformations in a way to reorder the elements on the main diagonal i.e. $g_r:=2+2\tilde{h}^2 \sn^2(k^{-1}r \tilde{h},k_{a,b})-(17/16) \tilde{h}^2$ in a non-decreasing order, without changing its eigenvalues. We will call $\hat{\ue{\ue{\mathfrak{T}}}}$ the transformed matrix and $\hat{g}_{r}$ the newly obtained sequence on the main diagonal. Let us suppose that $q_r>0$ (property that can be verified \emph{a posteriori}), in this way we can define $\gamma_r:=q_{r+1}/q_r$ and ensure that the new sequence $\gamma_{r+1}=\hat{g}_r-\gamma_r^{-1}$, with $\gamma_1=2-(17/16) \tilde{h}^2$, satisfies $\gamma_r \geq 1$ for all $r \leq \tilde{n}$. Such a sequence is commonly known as \emph{Sturm ratio sequence}. Knowingly, the latter cannot be studied by means of elementary methods due to the variation of $\hat{g}_r$. However, it is easy to show by induction that if we construct a sequence $\hat{g}_r' \leq \hat{g}_r$ then the sequence $\gamma_r'$ defined by $\gamma_{r+1}'=\hat{g}_r'-(\gamma_r')^{-1}$, $\gamma_1':=\gamma_1$ satisfies $\gamma_r'<\gamma_r$ for all $r \leq \tilde{n}$. \\
It can be checked that a suitable candidate for $g_r'$ is the following 
\beq{eq:ghatr}
\hat{g}_r':=
\begin{cases}
2-(17/16)\tilde{h}^2 & \qquad  1 \leq  r \leq \lceil 0.52 \tilde{h}^{-1}\rceil \\
2-(9/16)\tilde{h}^2 & \qquad \lceil 0.52 \tilde{h}^{-1}\rceil <  r \leq 
\lceil 1.95 \tilde{h}^{-1}\rceil\\
2+(\tilde{h}^2/16) &  \qquad
\lceil 1.95 \tilde{h}^{-1}\rceil < r  \leq \tilde{n}
\end{cases}
\mx{.}
\eeq
For instance, as for the first value, one needs to find $r^*$ such that $g_{r^*} \geq 2-(9/16)\tilde{h}^2$. By using the inequality $\sn (u,k) >\tanh u$, the latter is implied by $\tanh(k^{-1}r \tilde{h}) \geq 1/4$, which holds true (in particular) by choosing $r^*=0.26 \tilde{h}^{-1}$. Taking into account that each value of $\hat{g}_r$ appears twice due to the reordering, the first line of (\ref{eq:ghatr}) immediately follows. The second one is similar.\\
It is important to stress that a two-valued $\hat{g}_r'$ (i.e. the first and the last line of (\ref{eq:ghatr}) only) would be not sufficient for our purposes: basically, the approximation given by such a $\hat{g}_r'$ would be ``too rough'' for $\hat{g}_r$. A three-valued $\hat{g}_r'$ as in (\ref{eq:ghatr}) turns out to be the simplest acceptable option.\\
In order to ensure that $\gamma_r \geq 1$ it is sufficient to show that $\gamma_r' \geq 1$ for all $r \leq \lceil 1.95 \tilde{h}^{-1} \rceil $. This is due to the fact that the discrete dynamics associated to the third value of $\hat{g}_r'$, i.e. $\gamma_{r+1}'=2+(\tilde{h}^2/16)-(\gamma_r')^{-1}$ possesses an attractive fixed point strictly greater than one: in particular, any $\gamma_{\lceil 1.95 \tilde{h}^{-1}\rceil}' \geq 1$ is attracted by it, hence $\gamma_r' > 1$ for all $r > \lceil 1.95 \tilde{h}^{-1}\rceil$. \\
This implies that the analysis is reduced to the two recurrence equations  
\beq{eq:firstdyn}
\gamma_{r+1}'=\gamma_1'-(\gamma_{r}')^{-1}, \qquad 1 \leq r \leq \lceil 0.52 \tilde{h}^{-1}\rceil \mx{,}
\eeq
where $\gamma_1':=2-(17/16)\tilde{h}^2$, and, after an obvious redefinition for the purpose of shifting $r$, 
\beq{eq:seconddyn}
\begin{cases}
\gamma_{r+1}''& := 2-(9/16)\tilde{h}^2-(\gamma_{r}'')^{-1} \qquad 1 \leq r \leq\lceil 1.43 \tilde{h}^{-1}\rceil\\
\gamma_1'' & \leq \gamma_{\lceil 0.52 \tilde{h}^{-1}\rceil}'
\end{cases} 
\mx{.}
\eeq
This approach has a clear advantage, as the variation of $\hat{g}_r$ is removed, the solutions to (\ref{eq:firstdyn}) and (\ref{eq:seconddyn}) can be explicitly written. More in general, given $2\delta_1>g>0$ and the following recurrence equation
\[
\delta_{r+1} := g-(\delta_{r})^{-1} \mx{,}
\]
its solution reads as
\beq{eq:soldiscr}
\gamma_r=2^{-1}\left[g+  \sqrt{4-g^2} \cot (r \theta+\theta_1)) \right] \mx{,}
\eeq
where
\[
\theta:=\arctan  \left(g^{-1}\sqrt{4-g^2}\right), \qquad \theta_1:=\arctan \left((2\delta_1-g)^{-1}\sqrt{4-g^2} \right) \mx{.}
\]
As for (\ref{eq:firstdyn}), we have in particular $g=\gamma_1'$, hence $\theta_1=\theta$. By using the assumption on $\tilde{h}$ we easily get $2.06 \tilde{h} < \sqrt{4-(\gamma_1')^2}< 2.08 \tilde{h}$. As a consequence, we can bound $\tan[(j+1)\theta] \leq  \tan [1.04\tilde{h}(0.52\tilde{h}^{-1}+2) ]<0.63$. Hence, 
\beq{eq:gammafirststep}
\gamma_{\lceil 0.52 \tilde{h}^{-1}\rceil}' \geq 1-(17/32)\tilde{h}^2+(1.03/0.63)\tilde{h} > 1 + (1.63 - 0.54 \tilde{h})\tilde{h}=:\gamma_1'' \mx{.}
\eeq
At this point we can study the dynamic given by (\ref{eq:seconddyn}). In this case $g=2-(9/16)\tilde{h}^2$, so we get $1.31 \tilde{h} \leq \sqrt{4-g^2} \leq 1.5 \tilde{h}$. On the other hand $2 \gamma_1''-g=(13/4-\tilde{h}/2 )\tilde{h}>3 \tilde{h}$. This gives the required bound from below for the $\cot(\cdot)$ appearing in (\ref{eq:soldiscr}). In conclusion, 
\[
\gamma_{\lceil 1.43 \tilde{h}^{-1}\rceil}'' \geq 1-(9/32)\tilde{h}^2+(3/4)\tilde{h} \cot\{(1.43\tilde{h}^{-1}+1)\arctan[(1.5 \tilde{h})(2-9 \tilde{h}^2/16)]+\arctan(1/2)\}\geq 1
\]   
as it can be easily checked, in particular, for all $\tilde{h} \in (0, 0.01]$. 
\endproof
Proposition \ref{prop:defpos} plays a key role in the proof of the main result of this section, which states as follows
\begin{lem}\label{lem:bvp} Under the assumptions of Proposition \ref{prop:defpos}, define  
\[
Y_{r_0}:=\{(\ue{v},\ue{w})\in \mathfrak{D}(\ml{M}_{a,b})\times \mathfrak{D}(\ml{M}_{a,b})\, : \, ||\ue{v}||_{\infty},||\ue{w}||_{\infty} \leq r_0 \} \mx{,}
\] 
and set
\begin{align}
\mu_0 & :=\tilde{\alpha}^2 \pi^2 [8 \Delta_{a,b}^2(8 \Delta_{a,b}^2+3 \tilde{\alpha}^2)]^{-1} \mx{,} \label{eq:muzero} \\
r_0 & :=\pi^2 [8 \Delta_{a,b}^2+3 \tilde{\alpha}^2]^{-1} \label{eq:rhozero} \mx{.}
\end{align}
Then, for any given $\mu\in (0,\mu_0]$, there exists a unique $(\ue{v}^*,\ue{w}^*) \in Y_{r_0}$ satisfying (\ref{eq:discreteel}), which is the absolute minimum of the discretised functional $\ml{F}_{a,b}|_{\ml{M}_{a,b}}$.\\
Furthermore, for any arbitrarily chosen $(\ue{v}^{(0)},\ue{w}^{(0)}) \in Y_{r_0}$ the above mentioned solution satisfies $(\ue{v}^*,\ue{w}^*)=\lim_{k \rw \infty}(\ue{v}^{(k)},\ue{w}^{(k)})$, where $(\ue{v}^{(k)},\ue{w}^{(k)})$ is the sequence defined as follows 
\beq{eq:quasinewton}
(\ue{v}^{(k+1)},\ue{w}^{(k+1)})=(\underline{v}^{(k)},\underline{w}^{(k)})-\ue{\ue{\mathcal{J}_0}}^{-1} \underline{\Psi}(\underline{v}^{(k)},\underline{w}^{(k)},\ue{t};\mu), \qquad k=0,1,\ldots
\eeq 
\end{lem}
\begin{rem}\label{rem:remfinlemma}
As it is evident from (\ref{eq:muzero}), the threshold $\mu_0$  decreases as the time interval $\Delta_{a,b}$ increases, roughly, $\mu_0 \sim \Delta_{a,b}^{-4}$.  
\end{rem}
\proof The proof is based on a ``quantitative version'' of the classical IFT. The one contained in \cite{chier} will be used here. The conditions required by the latter, see \cite[(7.6),(7.7)]{chier}, read in our case as
\begin{align}
\sup_{\mu \in (0,\mu_0]} \norm{\mu \tilde{h}^2(\ue{V}(\ue{0},\ue{0},t),\ue{W}(\ue{0},\ue{0},t))} & \leq r_0 \left(2 \norm{\ue{\ue{\mathcal{J}_0}}^{-1}} \right)^{-1} \mx{,}\label{eq:settepsei}\\
\sup_{(\mu,\ue{v},\ue{w}) \in (0,\mu_0]\times Y_{r_0}}
\norm{\ue{\ue{\mathcal{I}}}-\ue{\ue{\mathcal{J}_0}}^{-1}\left[\ue{\ue{\ml{T}}}+\tilde{h}^2 (\ue{\ue{\ml{D}}}+\mu \ue{\ue{\ml{E}}} )\right]} & \leq 1/2
\label{eq:settepsette}
\mx{.}
\end{align}
The crucial point lies in finding an upper bound for 
$\norm{\ue{\ue{\mathcal{J}_0}}^{-1}}$. Despite this could turn out to be a very difficult task, at least in general, the previously stated positive definiteness of $\ue{\ue{\mathcal{J}_0}}$ and the fact that all its off-diagonal entries are negative, guarantee that $\ue{\ue{\mathcal{J}_0}}$ is an M matrix, see e.g. \cite[Pag. 626]{meyer2000matrix}. Hence, there exists $R \in \RR^+$ such that it can be written as
\beq{eq:mmat}
\ue{\ue{\mathcal{J}_0}}=R \ue{\ue{\mathcal{I}}} - \ue{\ue{\mathcal{B}}} \mx{,}
\eeq  
with $\ue{\ue{\mathcal{B}}}$ with non-negative entries and $\rho(\ue{\ue{\mathcal{B}}}) < R$. Note that this property can be recognized explicitly, simply by setting $R:=2+\tilde{h}^2$ and 
\[
\ue{\ue{\ml{B}}}:=
\begin{pmatrix}
	\underline{\underline{B}} & \ue{\ue{0}}\\
	\ue{\ue{0}} & \underline{\underline{B}}\\
\end{pmatrix}
+\tilde{h}^2
\begin{pmatrix}
	\underline{\underline{K}} & \ue{\ue{0}}\\
	\ue{\ue{0}} & \underline{\underline{0}}\\
\end{pmatrix} \mx{,}
\]
where $\ue{\ue{B}}:=\trid(1,0,1)$ and $\underline{\underline{K}}:=\diag(1-\cos q_0, 1-\cos q_1,\ldots)$. Clearly $\ue{\ue{\ml{B}}}$ has non-negative entries. On the other hand, as $\lambda^-(\ue{\ue{\mathcal{J}_0}})>(\tilde{\alpha}/\Delta_{a,b})^2 \tilde{h}^2 > 0$ by (\ref{eq:lowesteig}), it is immediate from (\ref{eq:mmat}) that $\rho (\ue{\ue{\ml{B}}})=R-\lambda^-(\ue{\ue{\mathcal{J}_0}})<R$ i.e. $\ue{\ue{\mathcal{J}_0}}$ is an M matrix. In particular we can bound, by using (\ref{eq:lowesteig}), 
\beq{eq:boundrho}
\rho (\ue{\ue{\ml{B}}}) \leq 2+ \Delta_{a,b}^{-2}(\Delta_{a,b}^2-\tilde{\alpha}^2) \tilde{h}^2=:\bar{\rho} (\ue{\ue{\ml{B}}}) \mx{.}
\eeq      
Furthermore, it is possible to show that, see e.g. \cite{meu}, the following inequality holds for all $k \geq 1$ and all $\ep>0$ 
\beq{eq:ineqb}
\norm{\ue{\ue{\mathcal{B}}}^k}^{\frac{1}{k}} \leq \rho (\ue{\ue{\ml{B}}}) + \ep \mx{.}
\eeq 
It is now sufficient to define, for instance, $\varepsilon :=(R-\bar{\rho} (\ue{\ue{\ml{B}}}))/2=(\tilde{\alpha}/\Delta_{a,b})^2 \tilde{h}^2/2$ and use (\ref{eq:mmat}) to write, see e.g. \cite[Pag. 618]{meyer2000matrix}, $\ue{\ue{\mathcal{J}_0}}^{-1}$ as a Neumann series: 
\begin{align*}
\norm{\ue{\ue{\mathcal{J}_0}}^{-1}} & \leq R^{-1} \sum_{k=0}^{+\infty} \norm{\left(R^{-1} \ue{\ue{\ml{B}}}\right)^k } \\
& \Heq{\ref{eq:ineqb}}{\leq} R^{-1} \sum_{k=0}^{+\infty} [R^{-1}(\rho (\ue{\ue{\ml{B}}}) + \ep)]^k\\
& \Heq{\ref{eq:boundrho}}{\leq} 2 (\Delta_{a,b}/\tilde{\alpha})^2 \tilde{h}^{-2} \mx{.}
\end{align*}
Now we can proceed with condition (\ref{eq:settepsei}). First of all let us note that, as can be easily deduced from (\ref{eq:jacobianpsi}), one has
\[
\norm{(\ue{V},\ue{W})}\leq 2 \mu_0 \tilde{h}^2 \mx{,}
\] 
hence, condition (\ref{eq:settepsei}) would require 
\beq{eq:condmuz}
\mu_0 \leq (\tilde{\alpha}/\Delta_{a,b})^2 r_0/8 \mx{,} 
\eeq
which holds true by virtue of (\ref{eq:muzero}) and (\ref{eq:rhozero}).\\ 
Let us now take into account condition (\ref{eq:settepsette}). Firstly, we expand 
$\ue{\ue{\ml{D}}}(\ue{v},\ue{t})=\ue{\ue{\ml{D}}}(\ue{0},\ue{t})+\ue{\ue{\ml{D}}}'(\ue{v},\ue{t})$ and observe that 
\beq{eq:consger}
\norm{\ue{\ue{\ml{D}}}'(\ue{v},\ue{t})+\mu\ue{\ue{\ml{V}}}(\ue{v},\ue{w},\ue{t})} \leq \max_{k=1,\ldots,\tilde{n}}\left(|v_{k}|+2\mu+\mu \right) \leq r_0+3 \mu_0 \Heq{\ref{eq:condmuz}),(\ref{eq:rhozero}}{\leq} (8 \Delta_{a,b})^{-2} \pi^2
\eeq
for all $(\ue{v},\ue{w}) \in Y_{r_0}$. Then we can proceed as follows
\begin{align*}
& \norm{\ue{\ue{\mathcal{I}}}-\ue{\ue{\mathcal{J}_0}}^{-1}\left( \ue{\ue{\mathcal{J}_0}} + \tilde{h}^2 \ue{\ue{\ml{D}}}'(\ue{v},\ue{t}) +\mu \tilde{h}^2 \ue{\ue{\ml{V}}}(\ue{v},\ue{w},\ue{t}) \right)}\\
 \leq &\tilde{h}^2 \norm{\ue{\ue{\mathcal{J}_0}}^{-1}} \norm{\ue{\ue{\ml{D}}}'(\ue{v},\ue{t})+\mu\ue{\ue{\ml{V}}}(\ue{v},\ue{w},\ue{t})}\\
 \leq & (2 \tilde{\alpha})^{-2} \pi^2 = 4/9 \mx{,}
\end{align*}
(recall the definition of $\tilde{\alpha}$ in the statement) hence (\ref{eq:settepsette}) is satisfied as $\Delta_{a,b} \geq 3 \pi$ and the thesis (with the exception of the minimality) follows from the IFT. \\
The proof is complete once the minimality of the solution has been shown. In order to achieve this step, we have to guarantee that the lowest (positive) eigenvalue of $\ue{\ue{\mathcal{J}_0}}$ is not changed ``too much'' by the contribution of $\ue{\ue{\ml{D}}}'+\mu\ue{\ue{\ml{V}}}$, which is controlled by $r_0$ and $\mu_0$. \\  
A standard argument of eigenvalue perturbation for symmetric matrices is used, see e.g. 
\cite[Theorem 8.1.5]{golub}. In particular we have that 
\beq{eq:perteigs}
\lambda^-(\ue{\ue{\mathcal{T}}}+\tilde{h}^2 \{ \ue{\ue{\mathcal{D}}}+\mu \ue{\ue{\mathcal{E}}} \}) \geq \lambda^-(\ue{\ue{\mathcal{J}_0}})-|\lambda^+[\tilde{h}^2(\ue{\ue{\ml{D}}}'+\mu\ue{\ue{\ml{V}}})]| \mx{.}
\eeq
The last term follows directly from the Gershgorin Theorem: more precisely, by denoting with $d_{i,j}'$ the elements of $\ue{\ue{\ml{D}}}'$ and with $v_{i,j}'$ those of $\ue{\ue{\ml{V}}}$, we have $S_i:=\tilde{h}^2 \sum_{i \neq j} |d_{i,j}'+\mu v_{i,j}'| \leq \mu \tilde{h}^2$. Hence, 
\[
\lambda^+[\tilde{h}^2(\ue{\ue{\ml{D}}}'+\mu\ue{\ue{\ml{V}}})] \leq \tilde{h}^2 \max_i (|d_{i,i}'|+\mu |v_{i,i}'| +S_i ) \leq \tilde{h}^2 (\norm{\ue{v}}+3 \mu) \Heq{\ref{eq:condmuz}),(\ref{eq:rhozero}}{\leq}  8 (\tilde{n}+1)^{-2} \pi^2 \mx{.}
\]  
In conclusion, from (\ref{eq:perteigs}),
\[
\lambda^-(\ue{\ue{\mathcal{T}}}+\tilde{h}^2 \{ \ue{\ue{\mathcal{D}}}+\mu \ue{\ue{\mathcal{E}}} \}) \geq (7/16) \pi^2 (\tilde{n}+1)^{-2} \mx{,}
\] 
which completes the proof.
\endproof

\clearpage 
\section{Drift sequence and time estimate}\label{sec:driftseq}
We now proceed with the very last step of the proof, which consists in the construction of a suitable sequence of frequencies $\{\omega_{i,i+1}\}$. The requirements for the latter have been anticipated in Sec. \ref{sec:convexity} and they are very well known in the literature: one of them is the ``closeness of two consecutive frequencies'' (\ref{eq:hypdue}), while the other one is related to the construction of the (double) sequence $(T_i,Q_i)$. In fact, as Lemma \ref{lem:convexity} holds in any set $\mathfrak{C}_{\cdot,\pi/4}$ we have to ensure that ``we can jump from any connected component of $\mathfrak{C}_{\cdot,\pi/4}$ to another'' with an arbitrarily long time, i.e. for any starting point $(\hat{T},\hat{Q}) \in \RR^2$ and any $\hat{\tau}>0$, there exists $t^*>\hat{\tau}$ such that  $\phi_{(1,\omega_{i,i+1})}^{t^*} (\hat{T},\hat{Q}) \in \mathfrak{C}_{\beta_{1,2}^*,\pi/4}$ for some $\beta_{1,2}^* \in \ZZ^2$, where we have denoted $\phi_{(1,\omega_{i,i+1})}^{t} (T,Q) :=(T+t,Q+\omega_{i,i+1} t)$. Such a feature, also known as ``ergodization property'', is systematically used in this context and it completes the scenario of what in a geometric language we would call $\lambda-$lemma: this is a fundamental ingredient of any geometric argument in Arnold's diffusion, see e.g. \cite{AIHPA_1994__60_1_1_0}. \\
Knowingly, the ergodization property relies on the numeric features of $\omega_{i,i+1}$, that is why the $\omega_{i,i+1}$ are usually assumed to be Diophantine. Without addressing the profound implications of ``ergodic nature'' which arise when dealing with such a class of numbers, we only recall the (remarkable) property for which the flow $\phi_{(1,\omega_{i,i+1})}^{t} (T,Q)$ falls in any arbitrarily small neighbourhood of any given point, at the price of taking sufficiently large $t$. Explicit bounds are known on such a time, see e.g. \cite{dumas1991}, \cite{bgw98} and \cite{Dumas1999}. However, as  $\diam_-\mathfrak{C}_{\cdot,\pi/4}$ is fixed and $O(1)$, $\omega_{i,i+1}$ in form of rational numbers are acceptable as well, provided that they are ``non-resonant up to a sufficiently high order'' \cite{MR1996776}. This feature, used in the latter work in its full generality, takes a simple and fully constructive form in our case. \\
In fact, it is sufficient to recall that the (periodic) flow on $\TT^2$ with frequency $(1,p/q)$, where $p,q\in \NN \setminus \{ 0 \}$ such that $\mcd(p,q)=1$, could be systematically avoiding sets $\mathcal{O} \subset \TT^2$ such that $\diam_-(\mathcal{O})<2 \pi/q$: this is exactly what we want to prevent from happening. For this purpose, as $\diam_-(\mathfrak{C}_{\cdot,\pi/4})= \pi/2$, it is sufficient to ``remove'' from the set of potential candidates for $\omega_{i,i+1}:=p^{(i)}/q^{(i)}$ all the frequencies in $[0,1]$ such that $q^{(i)} \leq 4$, i.e. the set Farey numbers up to order $4$, usually denoted with $F_4$, see \cite{ainsworth2012farey}. The elements of $F_4$ in increasing order will be denoted with $f_m$ (hence $f_0\equiv 0$ and $f_6\equiv 1$). \\
As a consequence, the ``non-resonant set'' we are looking for is simply 
\[
\mathcal{W}:= [0,1] \setminus F_4 \mx{.}
\]   
Given any $\omega_{a,b} \in \mathcal{W}$, let us define $\ep:=\dist(\omega_{a,b} ,f_{m^*})$, where $f_{m^*}$ is the closest\footnote{Clearly we shall set $\ep:=(f_1+f_2)/2$ if there exist $f_1, f_2 \in F_4$ such that $\dist(\omega_{a,b} ,f_1)=\dist(\omega_{a,b} ,f_2)$.} element of $F_4$ to $\omega^*$. As it is reasonable to expect, the ergodization time associated with $\omega_{a,b}$, $T_{\mx{erg}}(\omega_{a,b} )$, is singular as $\ep \rw 0$. Despite an upper bound requires a proof in the general case, as the one presented in \cite{MR1996776}, in our case it is straightforward to check that the following bound holds
\beq{eq:terg}
T_{\mx{erg}}(\omega_{a,b} ) \leq 2 \pi \ep^{-1}.
\eeq   
It is sufficient to recall Rem. \ref{rem:remfinlemma}, to realise that this bound raises a substantial problem:  the solution of the BVP (\ref{eq:bvpuno}) where $\omega_{a,b}=:(Q_b-Q_a)/(T_b-T_a)$, cannot be found via Lemma \ref{lem:bvp} if $\omega_{a,b} \in \ml{W}$ is ``too close'' to an element of $F_4$. The reason for this issue lies in hypothesis (\ref{eq:muzero}) which is not satisfied if $\Delta_{a,b} \sim T_{\mx{erg}}(\omega_{a,b} ) $ is ``too large''. More precisely, by comparing (\ref{eq:terg}) with (\ref{eq:muzero}), we find that Lemma \ref{lem:bvp} does not hold if 
\beq{eq:epzeroforb}
\ep \leq \ep_0:= 8 \pi \tilde{\alpha}^{-1}\left[ \sqrt{9 \tilde{\alpha}^2 + 4 \pi^2 / \mu } - 3 \right]^{-\frac{1}{2}} \mx{,}
\eeq
recalling that $\tilde{\alpha}=3 \pi/4$. Let us now define, 
\beq{eq:chi}
\chi_m:=(f_m+\ep_0,f_{m+1}-\ep_0) \mx{,}
\eeq
for all $m=0,\ldots,5$ such that $f_{m+1}-f_m \leq 2\ep_0$, otherwise we shall set $\chi_m:=\emptyset$. Hence, the BVPs for which Lemma \ref{lem:bvp} holds are those such that $\omega_{a,b}$ belongs to the set 
\beq{eq:wdmu}
\ml{W}_{d,\mu}:=\bigcup_{m=0}^5 \chi_m \mx{,}
\eeq
defining in this way the set mentioned in the main statement. Clearly, by (\ref{eq:epzeroforb}), $\meas(\ml{W}_{d,\mu}) \rw 1$ as $\mu \rw 0$, i.e. every $\chi_i$ turns out to be non-empty if $\mu$ is sufficiently small.  \\
In other terms the set $\ml{W}_{d,\mu}$ is the interval $[0,1]$ with some ``holes'' which have a $O(\sqrt[4]{\mu})$ diameter.  This size is actually ``quite large'' if compared to $C\mu$, but non-trivial drifting solutions can still be constructed in any non-empty $\chi_m$ as their diameter is $O(1)$, which is, independent of $\mu$. 
\begin{rem}%
The set $\ml{W}_{d,\mu}$ is uniquely defined once $\mu$ has been fixed. However, it is important to stress that the structure of $\ml{W}_d$ would have been more ``involved'' in the presence of $\mathfrak{C}_{\cdot,\delta}$ with $\delta$ smaller than $\pi/4$. More precisely, as $\delta$ decreases, one needs to increase the order $n$ of the Farey set which has to be removed from $[0,1]$. As a consequence, $\meas(\ml{N})$ would increase as well, and the resulting $\ml{W}_{d,\mu}$ would exhibit a higher number of connected components, say $M^*$, but with a sensibly smaller size. Still, if $\delta$ is independent of $\mu$, as in our case, this does not prevent the construction of non-trivial drifting solutions as $\meas([0,1]\setminus \{\cup_{m=0}^{M^*}\chi_i\}) \rw 0$ if $\mu \rw 0$.   
\end{rem}
As a result of the above described construction, every $\omega \in \mathcal{W}_{d,\mu}$ is ``ergodizing'' i.e. the flow on $\TT^2$ with frequency $(1,\omega)$ starting from any point on $\TT^2$ will reach any $\mathfrak{C}_{\cdot,\pi/4}$ in a finite time. In addition, the solution of the BVPs constructed on them exists by Lemma \ref{lem:bvp}. This implies that, given $\omega_I<\omega_F$ in any non-empty $\chi_m$, in order to construct the desired sequence $\{\omega_{i,i+1}\}$ we just need to ``fill the space'' between $\omega_I$ and $\omega_F$ in such a way (\ref{eq:hypdue}) is satisfied. This is easily achieved by setting 
\beq{eq:tranchain}
\omega_{i,i+1}=\omega_I+(N-2)^{-1}(\omega_F-\omega_I)(i-1), \qquad i=1,2,\ldots,N-1 \mx{,}
\eeq
where $N$ is chosen as in (\ref{eq:choicen}) (note that $N$ is even). We can finally proceed with the construction of the set $\mathfrak{C}_{\mathcal{I},\pi/4}$ in such a way the assumptions of Lemmas \ref{lem:convexity} and \ref{lem:bvp} hold true. Recalling Remark \ref{rem:convexity}, this is achieved by constructing increasing sequences $\{T_i^0,Q_i^0\}$ such that $\mathfrak{m}(T_i^0),\mathfrak{m}(Q_i^0)<\pi/4$, (\ref{eq:hypuno}) is satisfied and 
\beq{eq:tsuf}
T_{i+1}^0-T_i^0 \geq 3 \pi \mx{,}
\eeq
see Prop. \ref{prop:defpos}. The use of the superscript ``$0$'' is motivated by the fact that $\{T_i^0,Q_i^0\}$ will be used as the starting guess of the Nesterov sequence (\ref{eq:graddesc}), recall (\ref{eq:notationz}).\\ The use of the results of Sec. \ref{sec:bvp} on every subset $[T_i,T_{i+1}]\times [Q_i,Q_{i+1}]$ completes the constructive argument. \\
First of all let us notice that, 
\[
T_{i+1}^0-T_i^0 \equiv 2 k_{i,i+1} \K{k_{i,i+1}} \geq 2 k_0 \K{k_0} \geq \mathcal{T}^{-}(\mu) \geq  \mathcal{T}^{-}(\bar{\mu}) >10>3 \pi \mx{,}
\]
hence (\ref{eq:hypuno}) implies (\ref{eq:tsuf}).\\
Let us now assume (without loss of generality) $T_1^0=Q_1^0=0$. Trivially, $\mathfrak{m}(T_1^0)=\mathfrak{m}(Q_1^0)=0$ (recall the definition of $\mathfrak{m}$ given in Remark \ref{rem:convexity}), hence we can suppose to have determined $T_i^0,Q_i^0$ and proceed by induction. \\
Let us set 
\beq{eq:nistar}
n_i^*:=1+\left \lceil 1/6+(T_i^0+k_0 \K{k_0})/\pi\right \rceil \mx{,}
\eeq
then consider, for all $n \geq n_i^*$, the functions 
\begin{align}
T_{i+1}(n)&:=2 \pi n - T_i^0 \label{eq:tip1}\\
Q_{i+1}(n)&:=Q_i^0+\omega_{i,i+1}(T_{i+1}(n)-T_i^0) \label{eq:qip1} \mx{.}
\end{align} 
First of all, supposing $\mathfrak{m}(T_i)=0$, we have $\mathfrak{m}(T_{i+1}(n))=0<\pi/4$ by definition for all $n$. \\Given the properties of $\omega_{i,i+1}$, for all $Q_i^0$ there exist infinitely many $n\geq n_i^*$ such that $\mathfrak{m}(Q_{i+1}(n))<\pi/4$: let $n_i^0$ be the smallest of them. The inductive step is completed by setting  
\beq{eq:newtdef}
T_{i+1}^0:=T_{i+1}(n_i^0), \qquad Q_{i+1}^0:=Q_{i+1}(n_i^0) \mx{.}
\eeq
It is immediate to check that the choice (\ref{eq:nistar}) is such that (\ref{eq:tsuf}) holds true, in such a way the hypotheses of Lemma \ref{lem:bvp} are satisfied as well.\\
The estimate (\ref{eq:shadow}) can now be finally proven. Given any $1 \leq i \leq N-1$, let us consider the BVP 
(\ref{eq:bvps}). By (\ref{eq:omegai}) and Lagrange's Theorem, there exists $T_{i,i+1}^* \in (T_i,T_{i+1})$ such that 
\[
\dot{Q}(T_{i,i+1}^*)=\omega_{i,i+1} \mx{.}
\]
The maximum variation of $\dot{Q}(t)$ from this value can be bounded by recalling the basic approximation formula
\[
\ddot{Q}=\mu (1-\cos q) \sin Q = \mu (1-\cos q_0(t)) \sin Q_0(t)+O(\mu^2) \mx{,}
\]
where $q_0(t)$ and $Q_0(t)$ denote the zero-th order approximation in $\mu$ of $q(t)$ and $Q(t)$, respectively.  \\
By integrating and taking the absolute values one can find
\begin{align*}
\left|\dot{Q}(t) - \omega_{i,i+1}\right| & \leq \mu \int_{-\K{k}}^{\K{k}} \cn^2(\tau,k) d \tau \\
& = \mu k^{-2} \left[ \E{\am (\tau,k)}-(k')^2 \tau\right]_{-\K{k}}^{\K{k}}\\
& = 2 \mu k^{-2} \left( \E{k}-(k') \K{k} \right)\\
& \leq 2 \mu k^{-2} \E{1}\\
& \leq \mu \pi k^{-2}\\
& \leq \mu \pi + o(\mu) \mx{.}
\end{align*}
The last step consists in bounding how much the value of $\omega_{i,i+1}$ can vary during the minimisation process and this can be done via the following elementary geometrical argument: by definition $\omega_{i,i+1}$, is the slope of the line passing thorough the points $(T_i,Q_i)$ and $(T_{i+1},Q_{i+1})$. However, as we have shown by convexity, the points cannot move outside the boxes of size $2\delta$. In addition, the two boxes are, at least $\mathcal{T}^{-}(\mu)$ apart each other. It is easy to show that the maximum variation of the slope of all the possible lines is, say, $\omega_{max}-\omega_{min} \leq 4 \delta / \mathcal{T}^{-}(\mu)$. Clearly, as we solve the BVP in its discrete form, the bound has to take into account of the error arising from the discretisation, which is an $o(h)$ in our case. This proves (\ref{eq:shadow}).

\subsubsection*{Estimate of the time of drift along $\{\omega_{i,i+1}\}_{i=1,\ldots,N}$:} Let us choose a connected component $\chi_{m}\subset\mathcal{W}_{d,\mu}$  then, for all $\omega_I,\omega_F \in \chi_{n}$, with $\omega_I<\omega_F$. Set $\mathfrak{d}:=\omega_{F}-\omega_I$. The estimate of the time spent along the sequence $\{\omega_{i,i+1}\}$ exhibits a problem: as $\omega_I$ and $\omega_F$ could be placed anywhere in $\chi_m$, the distance from the closest Farey number does not depend continuously on $i$. A way to overcome this problem is to construct a fictitious sequence $\{\tilde{\omega}_{i}\}$ which represents the ``worst'' case of $\{\omega_{i,i+1}\}$ in terms of ergodization time. It is immediate to realise that such a sequence is symmetric with respect to the centre of $\chi_m$ and it can be defined as 
\[
\tilde{\omega}_i:=
\begin{cases}
\tilde{\omega}_I+\mathfrak{d}(2N-4)^{-1}(N-2i) &\qquad i=1,\ldots,N/2 \\
\tilde{\omega}_F-\mathfrak{d}(2N-4)^{-1}(N-2i) &\qquad i=N/2+1,\ldots,N
\end{cases} \mx{,}
\]   
where $\tilde{\omega}_I:=f_m+\ep_0$ and $\tilde{\omega}_F:=f_{m+1}-\ep_0$. \\
By (\ref{eq:terg}) a bound on the total ergodization time required is given by
\beq{eq:estimateone}
T_d(\{\tilde{\omega}_i\})=4 \pi \sum_{i=1}^{N/2}\frac{1}{\ep_0+\mathfrak{d}(2N-4)^{-1}(N-2i)} \leq 4 \pi (N-2)\int_0^{d/2}\frac{dx}{\ep_0+x} \mx{.}
\eeq
On the other hand, from (\ref{eq:epzeroforb}), one gets
\[
(2 \ep_0)^{-1} \leq \left(\sqrt{3}/64\right)\sqrt{9+8\mu^{-1/2}}\leq \left(8 \sqrt[4]{\mu} \right)^{-1} \mx{.}
\]
By using the latter, we find from (\ref{eq:estimateone}) 
\[
T_d(\{\tilde{\omega}_i\})\leq - \pi N \log (256 \mu) \mx{.}
\]
However, we have to recall that the time for a single transition is bounded by $\mathcal{T}^+(\mu)$. Hence, (\ref{eq:drifttime}) holds.

\begin{rem}\label{rem:delta}%
The importance of choosing $\delta$ ``large'' is now evident: a choice of a smaller $\delta$ would have required the use of a higher order of $F_n$. It is easy to realise that, given a segment $[\omega_I,\omega_F]$ the resonances to be crossed contained in this interval increase with $n$. In this way, each $\omega_{i,i+1}$ is closer to some resonance. A comparison with (\ref{eq:terg}), shows that the time $T_{\mbox{erg}}$ associated with each $\omega_{i,i+1}$ will be necessarily larger, and hence $T_{i+1}-T_i$.   
\end{rem}
\clearpage 
\section{Implementation and numerical experiments}\label{sec:experiments}
This section aims to present some of the most significant outcomes of the algorithm implementation. The latter has been outlined via a pseudo-code and reported in Appendix B. for the reader's convenience. Out of all the simulations performed, the analysis will be focused on three of them, in order to emphasise some key features of the constructed trajectories. The corresponding problem data are reported in Table \ref{tab:one}, whilst Table \ref{tab:two} collects some key information regarding the implementation performance, alongside a validation of bounds (\ref{eq:tmenopiu}).\\
As anticipated, an estimate of the computational cost of the whole task represents a key aspect as it could suggest whether or not the variational approach might be convenient with respect to computational methods of different nature. It is worth pointing out that the simulations of this work have been performed on a single-CPU personal machine, which means that the performance data in Tab. \ref{tab:two} have the potential to be remarkably improved. In particular, as suggested by Algorithm \ref{alg:two}, the $N$ BVPs are solved in a sequential fashion, whilst this could be efficiently dealt with, for instance, on machines with parallelisation capabilities.  \\  
Tab. \ref{tab:two} suggests that each Nesterov step (indexed by $r$) takes roughly $30 s$ every $10^3$ transitions, which leads to an estimate of about $30 ms$ in order to solve a single BVP. Given the ``high'' dimension of a single BVP (which has typically involved $n=O(10^4)$ variables) the performance time looks quite encouraging: as already mentioned in the introduction, the approximation of the BVP solution can rely on the very efficient ``quasi-Newton'' method (\ref{eq:quasinewton}) which can be used under the quantitative IFT hypotheses, see Lemma \ref{lem:bvp}.  
\begin{figure*}[h!]\begin{center}
		\vspace{0pt}
		\begin{minipage}[c][1\width]{0.49\textwidth}
			\hspace{0pt}
			{\begin{overpic}[width=\textwidth]{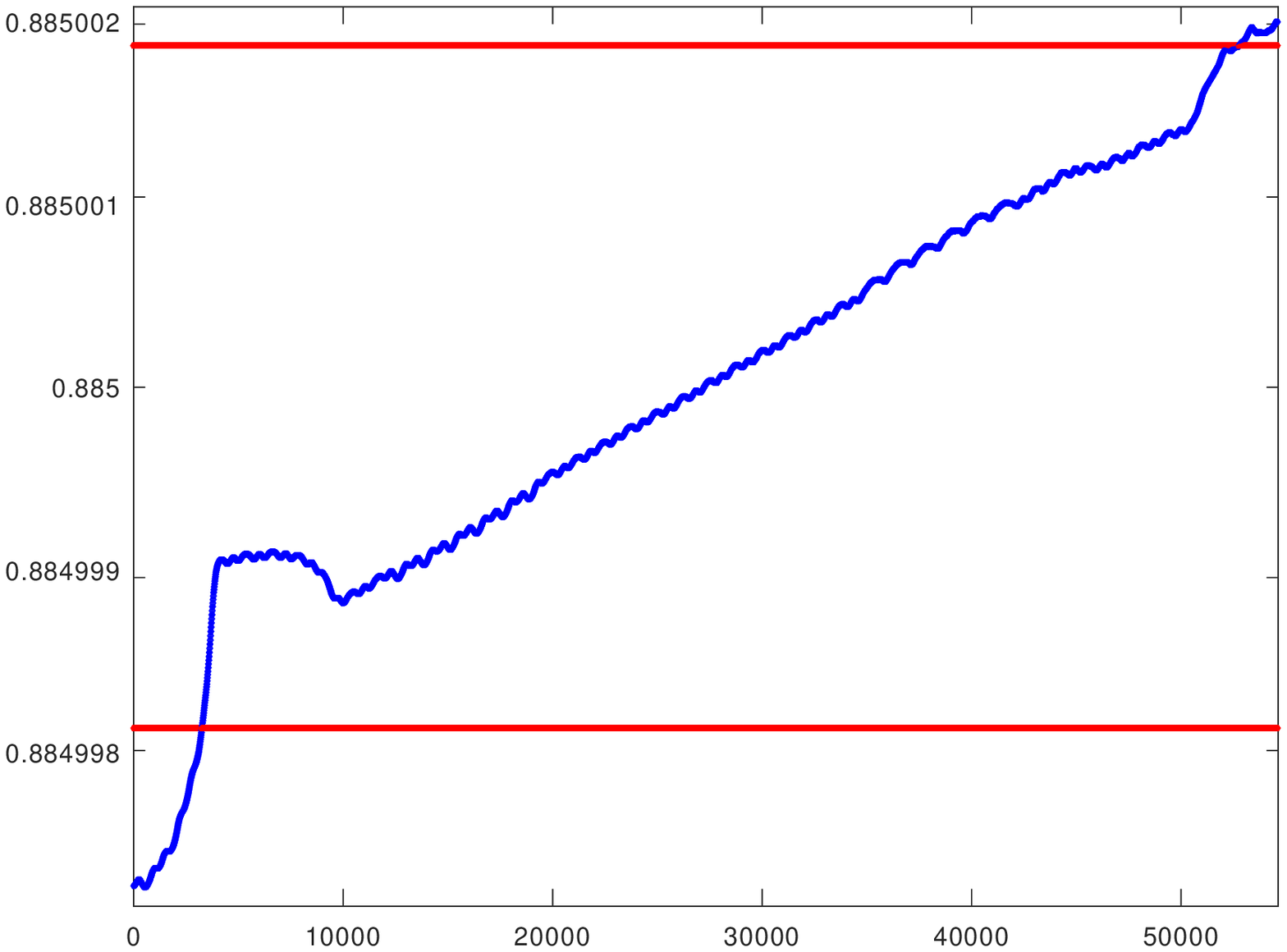}
					\put(00,75){\footnotesize (a)}
					\put(02,61){\footnotesize $\dot{Q}$}
					\put(75,2){\footnotesize $t$}
					\put(45,47){\footnotesize $\dot{Q}(t)$}
					\put(66,16){\footnotesize $\omega_I$}
					\put(66,63){\footnotesize $\omega_F$}
			\end{overpic}}
		\end{minipage}	
		\begin{minipage}[c][1\width]{0.49\textwidth}
			\hspace{0pt} 
			{\begin{overpic}[width=\textwidth]{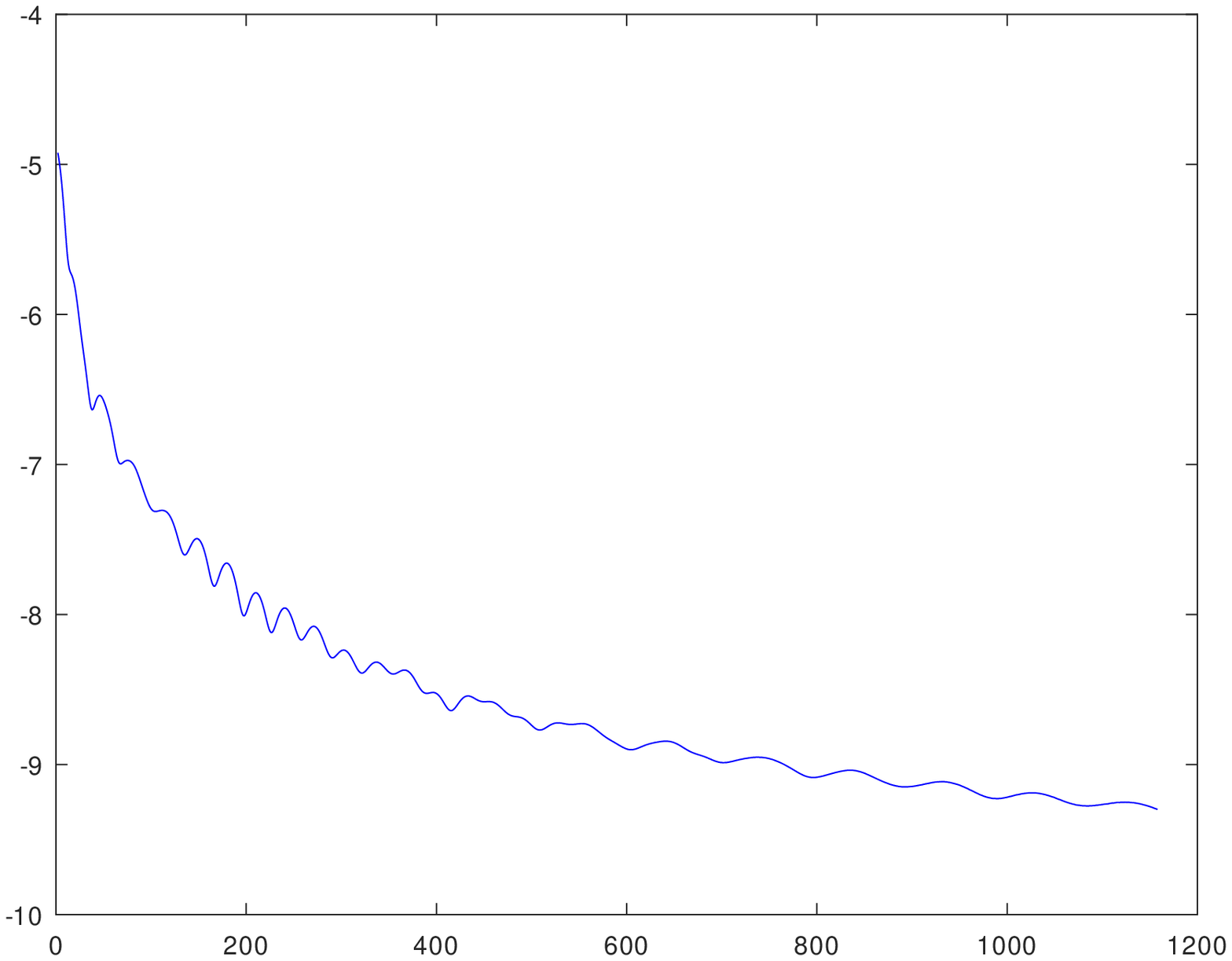}
					\put(00,75){\footnotesize (b)}
					\put(75,2){\footnotesize $r$}
					\put(40,28){\footnotesize $\log_{10}|\nabla \ml{F}|$}
			\end{overpic}}
		\end{minipage}\\
		\vspace{-20pt}
		\caption{Plot of the drifting trajectory of Simulation $n.1$, depicted in panel (a), and behaviour of the corresponding $|\nabla \ml{F}(Z_r)|$ in a logarithmic scale during the minimisation dynamics (panel (b)). It is worth emphasising the behaviour of $\dot{Q}$ in the time interval $t \in [0, \sim 10^4]$, characterised by a faster drift (w.r.t. the rest of the dynamics), followed by a stage in which the dynamics is ``nearly-stationary'' and finally, which is the most interesting aspect, a stage in which it even decreases. Clearly, 
$\dot{Q}$ is here the shorthand for $\dot{\ue{Q}}^{R,(1,N+1)}$, where $R$ is the number of steps performed by the Nesterov method.}
		\label{fig:one}
		\vspace{0pt}
	\end{center}
\end{figure*}
\\
In all the three cases reported, she perturbation size has been set to $\mu=0.75\cdot 10^{-7}$. This is arguably a ``very small'' value. However, as it can be easily verified with an explicit calculation for the cases at hand, the threshold (\ref{eq:muzero}) only allows values which are just slightly over $10^{-7}$. For instance, in Sim. n. $3$, this value is set to $1.34 \cdot 10^{-7}$. The reason beyond this restriction lies in the necessity to solve the BVP over ``long'' time intervals, see also Rem. \ref{rem:remfinlemma} , requiring in this way ``very small'' perturbation sizes, in order to guarantee the convergence of the method described in Lemma \ref{lem:bvp}.\\ 
Some comments are in order concerning the properties of the constructed solutions. As it is clear from all the three plots, $\dot{Q}(t)$ is not monotonically increasing. Whilst this feature is evident from Fig.\ref{fig:one}(a), it can be witnessed in Fig.\ref{fig:two}(a),(b) as well, despite the motion clearly takes place in a region of the phase space in which it looks much more ``regular''. In fact, even in this case, a more closer look, shows a sequence of increasing-decreasing patterns for $\dot{Q}(t)$. More importantly, this behaviour does not repeat at each transition (as it would be much more natural to expect) but after ``several'' transitions. For instance, in the solution of Fig.\ref{fig:two}(a) it is possible to count slightly more than $60$ ``crests''. By recalling that we are dealing with $1002$ transitions in this case, such a behaviour repeats itself roughly every $17$ of them. 
\begin{figure*}[h!]\begin{center}
		\vspace{0pt}
		\begin{minipage}[c][1\width]{0.49\textwidth}
			\hspace{0pt}
			{\begin{overpic}[width=\textwidth]{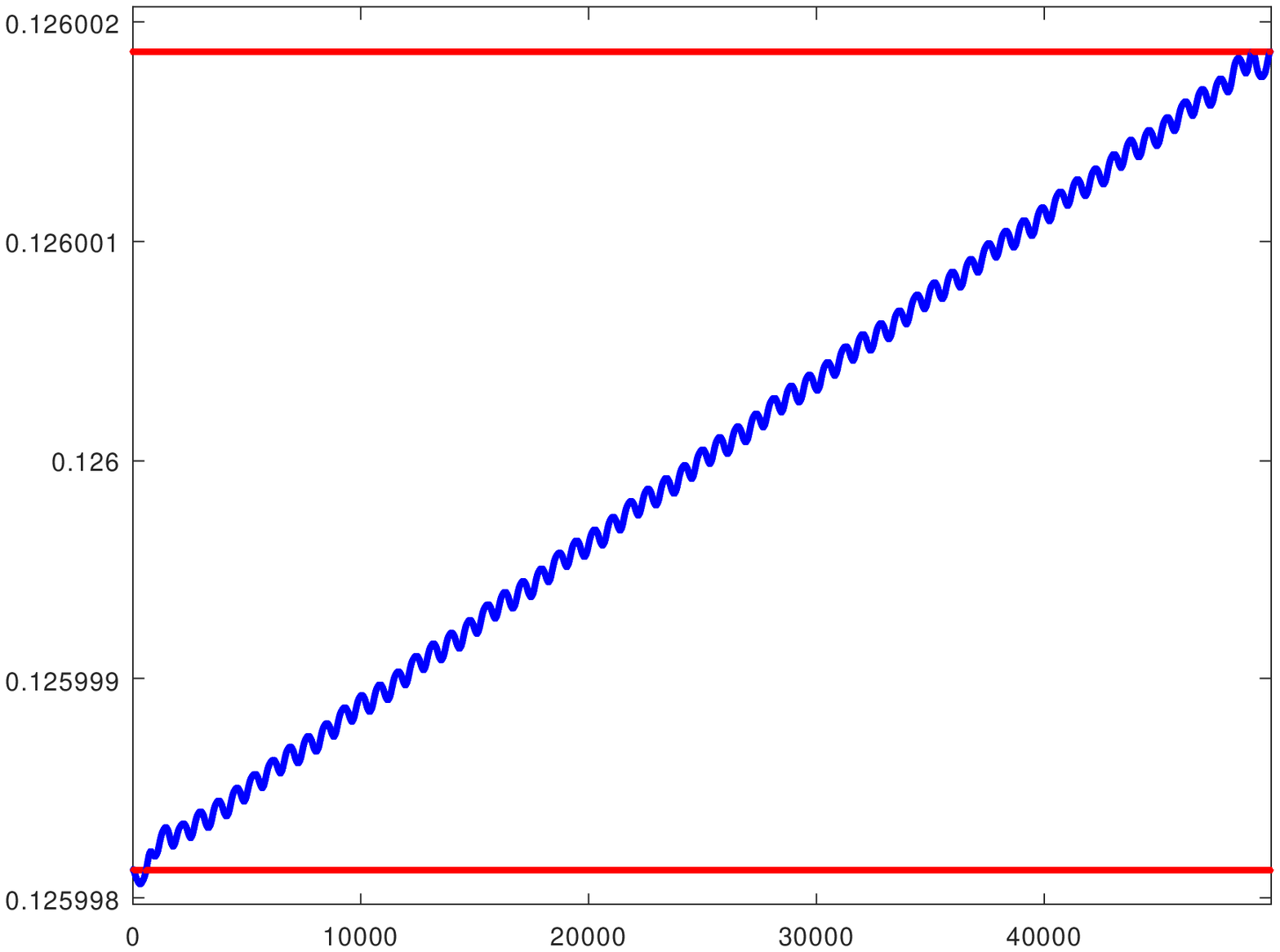}
					\put(00,75){\footnotesize (a)}
					\put(02,61){\footnotesize $\dot{Q}$}
					\put(75,2){\footnotesize $t$}
					\put(45,42){\footnotesize $\dot{Q}(t)$}
					\put(66,12){\footnotesize $\omega_I$}
					\put(66,63){\footnotesize $\omega_F$}
			\end{overpic}}
		\end{minipage}	
		\begin{minipage}[c][1\width]{0.49\textwidth}
			\hspace{0pt} 
			{\begin{overpic}[width=\textwidth]{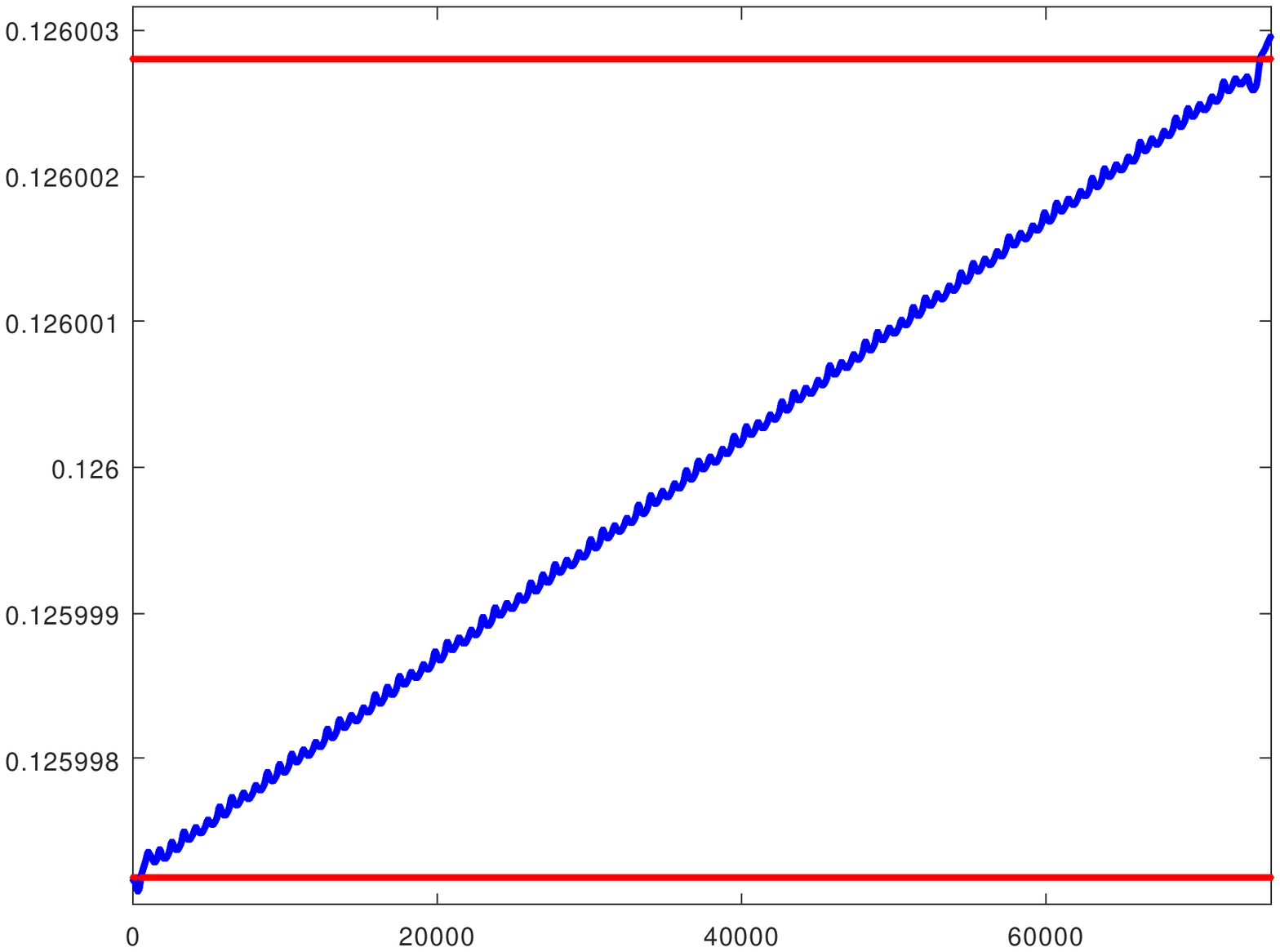}
					\put(00,75){\footnotesize (b)}
					\put(75,2){\footnotesize $t$}
					\put(45,42){\footnotesize $\dot{Q}(t)$}
					\put(66,12){\footnotesize $\omega_I$}
					\put(66,63){\footnotesize $\omega_F$}
			\end{overpic}}
		\end{minipage}\\
		\vspace{-20pt}
		\caption{The trajectories constructed in the case of Sim. n. $2$ and $3$, respectively. Despite the drift looks ``much more regular'' with respect to the case depicted in Fig. \ref{fig:one} (a), a closer look shows a sequence of ``crests''.}
		\label{fig:two}
		\vspace{0pt}
	\end{center}
\end{figure*}
\\
This ``microscopic'' behaviour seems to appear, on a ``macroscopic'' scale, in the first part of the trajectory depicted in Fig. \ref{fig:one}(a), roughly, in the time interval $t \in [0, \sim 10^4]$. The behaviour looks similar: a faster (w.r.t. the average slope) ``increase'' is followed by a ``stationary'' and a final ``decreasing'' stage. Interestingly enough, this would rather resemble a stochastic behaviour as those investigated, e.g., in the already cited papers \cite{MR2004513},\cite{efthy}. 
This represents a substantial difference with respect to the early geometrical methods (the very same original Arnold's mechanisms included), in which the dynamic increases the value of $\dot{Q}(t)$ \ue{at each step}, in a monotonic fashion. The latter approach has been beautifully described in \cite{10.1007/978-1-4020-6964-2_13} with a ``sailing'' metaphor:
\begin{quotation}
``When the wind is favorable, the boat moves. When the wind turns bad, it moves to the coast and anchors waiting till the wind becomes favorable again.''
\end{quotation}
In the light of the simulations reported, one could try to rephrase this for the ``variational sailor'', as the one who has the ability to (``plan'' how to) take advantage of favorable as well as unfavorable winds (and possibly anchoring) in the ``best'' way to reach their destination.\\
\bigskip 
\bigskip
\begin{table}[H]
\begin{center}
\begin{tabular}{|c|l|l|l|l|l|}
\hline
	Sim. n. & $\omega_I$ & $\omega_F$ & $\lceil \mu^{-1}(\omega_F-\omega_I) \rceil $ & N. trans. &  $T_d/N$ \\
\hline
	1 & 0.884998 & 0.885002 & 50 & 1000 & 54.632 \\
\hline
	2 & 1.259981 &  1.260019 & 50 & 1002 & 49.864 \\
\hline
	3 & 0.125997 & 1.260028 & 75 & 1500 & 49.867\\
\hline
\end{tabular}
\captionof{table}{Main problem data of the simulations reported. In the last column, the effective time of drift per transition. The latter can be compared with the bounds (\ref{eq:tmenopiu}) which, for the chosen value of $\mu=0.75 \cdot 10^{-7}$, give $\ml{T}^{-}(\mu)=16.631$ and $\ml{T}^{+}(\mu)=71.840$.   \label{tab:one}}
\end{center}
\end{table}

\begin{table}[H]
\begin{center}
\begin{tabular}{|c|l|l|l|l|l|}
\hline
	Sim. n.  & Nesterov steps & Elaps. t. (s) & Elaps. t. per BVP (ms)\\
\hline
	1  & 1159 &  30135.2 & 26.001\\
\hline
	2  & 166 & 5215.99 & 31.359\\
\hline
	3  & 134 & 6543.5 & 32.555\\
\hline
\end{tabular}
\captionof{table}{Some performance data of the three simulations. In particular, in the third column, the time required by the Nesterov algorithm to approximate the minimum (within a prefixed precision) is reported.  \label{tab:two}}
\end{center}
\end{table}

\section{Conclusions}
Via a suitable adaptation of proofs from existing works and complementing them with the necessary tools and ingredients, the paper shows how the variational approach introduced by U. Bessi and thereafter developed by himself and other authors, can be formulated on a discrete space and being suitable for a full machine implementation. Bessi's original idea to join together solutions of BVPs with suitable boundary conditions has been naturally implemented via a two-layer algorithm. This has led to several advantages in terms of CPU time and memory usage, primary concerns in large scale optimisation problems as the one obtained here. Performance data have been discussed, and non trivial problems with thousands of transitions have been shown to be compatible even with personal machines capabilities. \\ 
A special attention is deserved, however, by the behaviour of the constructed solutions. As shown in Sec. \ref{sec:experiments}., in fact, it appears far from being trivial: as it looks, the method is ``free to arrange'' the dynamics in increasing or decreasing patterns, like some sort of a ``suitably organised stochasticity''. On the other hand, diffusion as a form of instability, only requires a ``coarse'' shadowing of the transition chain of tori, or, possibly, no shadowing at all, provided that the initial and final values of $\dot{Q}$ are separated as prescribed (we remark that, for instance, \cite{MR1996776}, uses $\omega_{I,F}$ as boundary conditions of a suitable functional). As such, the variational approach seems to be particularly adequate in taking advantage of this freedom. 

\clearpage 
\section*{Appendix A. Pendulum and elliptic functions}\label{app:a}
Below, a short reminder on the application of Jacobian elliptic functions to the simple pendulum is reported.  (Much) more details can be found, e.g., in  \cite{Brizard_2009}.   \\
Let us consider the following Lagrangian 
\beq{eq:pendlag}
L=\frac{\dot{q}^2}{2}+(1-\cos q) \mx{,}
\eeq
corresponding to a simple pendulum with unstable equilibria located at $2j\pi$, $j \in \ZZ$. Let us notice that this feature makes it different from other models usually found in the literature. The argument recalled below is clearly analogous but it gives rise to slightly different expressions for (\ref{eq:am}) and (\ref{eq:kappa}).\\
For instance, the rotating solutions of (\ref{eq:pendlag}) can be written as 
\beq{eq:am}
q(t)=\pi+2 \am(k^{-1}t,k) \mx{,}
\eeq
where $\am(u,k)$ is the \emph{Jacobi amplitude} and $k$ is a parameter known as \emph{modulus}. We recall that, for any $k \in [0,1]$ and $u \geq 0$, the ``amplitude of $u$'' is defined as the value of the angle $\varphi$ such that $u=F(\varphi,k)$, where 
\[
F(\varphi,k)=\int_0^{\varphi} \frac{dx}{\sqrt{1-k^2 \sin^2 x}} \mx{,}
\]
and it is written $\varphi=\am(u,k)$. The latter give rise to the definition of the so-called \emph{copolar trio}
\[
\sn (u,k):=\sin \left(\am(u,k)\right), \qquad \cn(u,k):=\cos \left(\am(u,k)\right),\qquad \dn(u,k):=\sqrt{1-k^2 \sn^2(u,k)} \mx{.}
\]
As for the model (\ref{eq:pendlag}), the parameter $k$ is related to the energy $2 h:=\dot{q}^2/2+(\cos q -1)$ of the corresponding motion, via the relation
\beq{eq:kappa} 
k=1/\sqrt{1+h} \mx{.}
\eeq
Knowingly, rotations take place for all $k \in (0,1)$ (i.e. $h>0$) while the limit value $k=1$ (i.e. $h=0$), corresponds to the separatrix.  \\
By (\ref{eq:am}), the Lagrangian (\ref{eq:pendlag}) reads as
\beq{eq:lagraell}
L=2 k^{-2} \dn^2 (k^{-1}t,k)+\cn^2 (k^{-1}t,k) \mx{,}
\eeq
while the conservation of energy, via (\ref{eq:kappa}), is reduced nothing but to the Jacobi identities $\sn^2(u,k)+\cn^2(u,k)=1$ and $k^2 \sn^2(u,k)+\dn^2(u,k)=1$.
\\
The period of the rotation is expressed in term of the \emph{complete elliptic integral of the first kind} 
\beq{eq:period}
T:=2 k \K{k}, \qquad \K{k}:=F(\pi/2,k) \mx{.}
\eeq
The latter, defined with respect to the \emph{complementary modulus} $k':=\sqrt{1-k^2}$, is usually denoted as $\Kp{k}:=\K{k'}$. Finally, let us recall the \emph{complete elliptic integral of the second kind}, which is defined as 
\[
\E{k}:=\int_0^{\frac{\pi}{2}} dx \sqrt{1-k^2 \sin^2 x} \mx{.}
\]

\section*{Appendix B.}

\begin{algorithm}[H]%
\caption{Construction of drifting trajectories}
\label{alg:one}
\begin{algorithmic}
\State \textbf{Input:} $\mu>0$ ``sufficiently small'' and $\epsilon>0$ (tolerance);  
\State $\bullet$ Compute the set $\ml{W}_{d,\mu}$ via (\ref{eq:wdmu}), (\ref{eq:chi}) and (\ref{eq:epzeroforb}); 
\State $\bullet$ Pick $\omega_I<\omega_F$ in the same connected component of $\ml{W}_{d,\mu}$; 
\State $\bullet$ Compute the number $N$ of transitions needed according to (\ref{eq:choicen}) ;
\State $\bullet$ Compute the sequence $\{\omega_{i,i+1}\}_{i=1,\ldots,N-1}$ ``transition chain'' by using (\ref{eq:tranchain}) ; 
\State $\bullet$ Compute the sequences $\{T_i^0,Q_i^0\}$ following the argument described in Sec., and in particular steps (\ref{eq:nistar}) $\rw$ (\ref{eq:newtdef}). The total time of drift is given by $T_{N+1}$;
\State $\bullet$ $\Delta_{a,b} \gets \max_{i=1,\ldots,N} T_{i+1}-T_i$ and compute $\mu_0$ via (\ref{eq:muzero}) ;
\If{$\mu \leq \mu_0$} proceed 
\EndIf  
\State $\bullet$ $Z_0 \gets (T_2^0,Q_2^0,\ldots,T_N^0,Q_N^0)$ i.e. (\ref{eq:notationz}) and $r \gets 0$, $\mathfrak{N}_{\ml{F}} \gets 2 \epsilon$; 
\While{$\mathfrak{N}_{\ml{F}}>\epsilon$}
\State $\circ$ Evaluate $\nabla \ml{F}(Z_r)$ via \textbf{Algorithm 2} then $\mathfrak{N}_{\ml{F}} \gets |\nabla \ml{F}(Z_r)|$;
\State $\circ$ Compute $Z_{r+1}$ from $Z_r$ via (\ref{eq:graddesc}), (Nesterov's algorithm);
\State $\circ$ $r \gets r+1$;
\EndWhile
\State $\bullet$ $r^* \gets r-1$; 
\State \textbf{Output:} $(T_2^*,Q_2^*,\ldots,T_N^*,Q_N^*) \gets Z_{r^*}$ is the (approximated) minimum of $\ml{F}$;   
\end{algorithmic}
\end{algorithm}

\begin{algorithm}[H]
\caption{Evaluation of $\nabla \ml{F}(Z_r)$}
\label{alg:two}
\begin{algorithmic}
\State \textbf{Input:} $Z_r$ and all the other parameters involved
\For{i=1,\ldots,N}
\State $\circ$ $T_a \gets T_i$, $T_b \gets T_{i+1}$ and define the other necessary objects as in Sec. 6;
\State $\circ$ Construct $(\ue{q}^0,\ue{Q}^0)$ as (discrete) solution of the unperturbed problem on $[T_i,T_{i+1}]$;
\State $\circ$ $(\ue{v}^{(0)},\ue{w}^{(0)}) \gets (\ue{0},\ue{0}) $, $k \gets 0$ and $\mathfrak{N}_{\ml{N}} \gets 2 \epsilon$;
\While{$\mathfrak{N}_{\ml{N}}>\epsilon$}
\State $\cdot$ Compute $(\ue{v}^{(k+1)},\ue{w}^{(k+1)})$ from $(\ue{v}^{(k)},\ue{w}^{(k)})$ via (\ref{eq:quasinewton}) (Quasi-Newton method);
\State $\cdot$ $\mathfrak{N}_{\ml{N}} \gets |\ue{v}^{(k+1)}-\ue{v}^{(k)}|+|\ue{w}^{(k+1)}-\ue{w}^{(k)}|$;
\State $\cdot$ $k \gets k+1$;
\EndWhile
\EndFor
\State $\bullet$ $k^* \gets k-1$; 
\State $\bullet$ $(\ue{q},\ue{Q}) \gets (\ue{q}^0+\ue{v}^{(k^*)},\ue{Q}^0+\ue{w}^{(k^*)})$; 
\State $\bullet$ Compute $(\dot{\ue{q}},\dot{\ue{Q}})$ in the discrete sense (finite differences); 
\State \textbf{Output:} Compute $\nabla \ml{F}(Z_r)$ via (\ref{eq:gradbertiuno}) and (\ref{eq:gradbertidue});
\end{algorithmic}
\end{algorithm}

\clearpage 
\subsection*{Acknowledgements}
I am extremely grateful to Prof. G. Gallavotti for having introduced me to the ``constructivity problem'' at the time I was one of his Ph.D. students, as well as to Prof. U. Bessi for several discussions about his cited works and remarkable comments.\\
I wish to thank those colleagues who have kindly shown interest in this problem and, in particular, to Prof. L. Biasco and Prof. A. Sorrentino.\\  
This research has been conducted during my affiliation with the School of Mathematics of the University of Bristol, UK, under the support of ONR Grant No. N00014-01-1-0769. I wish to thank Prof. S. Wiggins for stimulating discussions. However, some minor parts of the paper were unfinished and they have been completed only recently, whilst affiliated with the University of Naples Federico II, under the support from the Italian Ministry of University and Research, PRIN2020 funding program, grant number 2020PY8KTC. I am grateful to Prof. M. d'Aquino for this opportunity.\\ 
The numerical simulations reported in this paper and the corresponding plots have been performed with GNU Octave \cite{oct}.

\bibliographystyle{alpha}
\bibliography{Diff.bib}

\end{document}